\documentclass[article]{amsart}
\usepackage[T1]{fontenc}
\usepackage{amsthm}
\usepackage{amssymb}
\usepackage{amsmath}
\usepackage{mathabx}

\usepackage{xcolor}
\usepackage{verbatim, color}

\usepackage{bbm}
\usepackage{float}
\usepackage{dsfont, graphicx}
\usepackage{epstopdf}
\usepackage[hidelinks]{hyperref}
\usepackage[T1]{fontenc}
\usepackage{amssymb,amsmath, amsthm, amsfonts}
\usepackage{graphicx}
\usepackage{tikz}
\usepackage{tkz-euclide}

\usepackage[hidelinks]{hyperref}
\usepackage{pagerange}
\usepackage{thmtools,thm-restate}

\makeatletter
\numberwithin{equation}{section}
\numberwithin{figure}{section}
\theoremstyle{plain}
\newtheorem{thm}{Theorem}[section]
\newtheorem{prop}[thm]{Proposition}

\newtheorem{lem}[thm]{Lemma}

\newtheorem{rem}[thm]{Remark}

  \newcounter{casectr}
  
\theoremstyle{definition}

\theoremstyle{remark}

\newcommand{\R}{\mathbb{R}}
\newcommand{\T}{\mathbb{T}}

\newcommand*{\dif}{\mathop{}\!\mathrm{d}}

\def \bR {\Bbb R}

\def \bZ {\Bbb Z}

\def \bT {\Bbb T}

\def \l {\left}
\def \r {\right}

\def \lea {\lesssim}
\def \gea {\gtrsim}

\def \p {^\prime}

\DeclareMathOperator{\supp}{\mathrm{supp}}

\begin{document}

\address{Yangkendi Deng
\newline \indent Academy of Mathematics and Systems Science, Chinese Academy of Sciences, Beijing, China.}
\email{dengyangkendi@amss.ac.cn}

\address{Chenjie Fan
\newline \indent Academy of Mathematics and Systems Science and Hua Loo-Keng Key Laboratory of Mathematics, Chinese 
\newline \indent Academy of Sciences, Beijing, China.}
\email{fancj@amss.ac.cn}

\address{Kailong Yang
\newline \indent Chongqing National Center for Applied Mathematics, Chongqing Normal University,
  Chongqing, China. }
\email{ykailong@mail.ustc.edu.cn}

\address{Zehua Zhao
\newline \indent Department of Mathematics and Statistics, Beijing Institute of Technology, Beijing, China.
\newline \indent Key Laboratory of Algebraic Lie Theory and Analysis of Ministry of Education, Beijing, China.}
\email{zzh@bit.edu.cn}

\address{Jiqiang Zheng
\newline \indent Institute of Applied Physics and Computational Mathematics,
Beijing, 100088, China.
\newline\indent
National Key Laboratory of Computational Physics, Beijing 100088, China}
\email{zheng\_jiqiang@iapcm.ac.cn, zhengjiqiang@gmail.com}

\title{On bilinear Strichartz estimates on waveguides with applications}
\author{Yangkendi Deng, Chenjie Fan, Kailong Yang, Zehua Zhao and Jiqiang Zheng}
\maketitle

\begin{abstract}
We study local-in-time and global-in-time bilinear Strichartz estimates for the Schr\"odinger equation  on waveguides.  As applications,  we apply those estimates to study global well-posedness of nonlinear Schr\"odinger equations on these waveguides.
\end{abstract}

\tableofcontents
\section{Introduction}
\subsection{Statement of main results}
In this article, we study bilinear Strichartz estimates in the waveguide setting. Consider linear Schr\"odinger equations on rescaled waveguide $\mathbb{R}^{m}\times \mathbb{T}_{\lambda}^{n}$,
\begin{equation}\label{eq: lsmain}
\begin{cases}
iu_{t}+\Delta_{x,y} u=0,\\
u(x,y,0)=f(x,y), (x,y)\in \mathbb{R}^{m}\times \mathbb{T}_{\lambda}^{n}.
\end{cases}	
\end{equation}
Here $\mathbb{T}_{\lambda}=\lambda \mathbb{T}$ is the rescaled tori.

We denote the linear propagator of \eqref{eq: lsmain} via $U_{\lambda}(t)$, i.e. the solution $u$ to \eqref{eq: lsmain} satisfying $u=U_{\lambda}(t)f$. 

Our main result is
\begin{thm}\label{thm: bilinearmain}
Consider \eqref{eq: lsmain}. When $m=1, n=1$, one has
\begin{equation}\label{eq: bilinear1}	
 \| U_{\lambda}(t) P_{N_1}f  U_{\lambda}(t) P_{N_2}g \|_{L_{x,y,t}^2(\mathbb{R}\times \mathbb{T}_{\lambda} \times [0,1])} \lesssim \left(\frac{1}{\lambda} +\frac{N_2}{N_1}  \right)^{\frac{1}{2}} \|f\|_{L^2} \|g\|_{L^2}.
\end{equation}
When $m\geq 2, n\geq 1$, let $d=m+n$ be the full dimension, one has
\begin{equation}\label{eq: bilinear2}
	 \| U_{\lambda}(t) P_{N_1}f  U_{\lambda}(t) P_{N_2}g \|_{L_{x,y,t}^2(\mathbb{R}^{m}\times \mathbb{T}^{n}_{\lambda} \times [0,\infty))} \lesssim \left( \frac{N_2^{d-3}}{\lambda} +\frac{N_2^{d-1}}{N_1}  \right)^{\frac{1}{2}} \|f\|_{L_{x,y}^2} \|g\|_{L_{x,y}^2}.
\end{equation}
Here $N_{1}\ge N_{2}\geq 1$ are dyadic numbers, and $P_{N}$ is the Littlewood-Paley projection to frequency $\sim N$.
\end{thm}

\begin{rem}
It remains an interesting problem to further study the case $m=1, n\geq 2$.	
\end{rem}
 \begin{rem}
 Estimates \eqref{eq: bilinear1} and \eqref{eq: bilinear2} are sharp in the sense that one can construct examples saturating those estimates and \eqref{eq: bilinear1} can only hold locally even for $N_{1}\gg N_{2}$, we will present the examples in Appendix \ref{app: sharp}.
 \end{rem}

 As applications of \eqref{eq: bilinear1}, one may follow \cite{DPST} and use I-method of \cite{I-team2} to obtain low regularity GWP for cubic NLS on waveguide $\mathbb{R}\times \mathbb{T}$,
\begin{equation}\label{eq: cubicnls}
\begin{cases}
iu_{t}+\Delta_{x,y} u=|u|^{2}u,\\
u(x,y,0)=f(x,y)\in H^{s}(\mathbb{R}\times \mathbb{T})
\end{cases}
\end{equation}
 for $s>\frac{2}{3}$.

 We take this chance to improve a bit the computation in \cite{DPST}, and prove
 \begin{thm}\label{thm: nonlinearmain}
 Equation \eqref{eq: cubicnls} is {globally}  well-posedness for initial data in $H^{s}(\mathbb{R}\times \mathbb{T})$ with $s>\frac{3}{5}$.	
 \end{thm}
 \begin{rem}
 Since I-method does not capture $\epsilon$ loss, (because it is a subcritical method), our results in Theorem \ref{thm: nonlinearmain} also holds for cubic NLS on $\mathbb{T}^{2}$.	It is worth noting that Schippa \cite{schippa2023improved} recently achieved the same low regularity global well-posedness as ours for the defocusing cubic NLS on a two-dimensional torus by modifying the I-energy through the resonance decomposition.
 \end{rem}

We also further study $(2k+1)$-th order NLS on $\mathbb{R}\times \mathbb{T}$
\begin{equation}\label{eq: knls}
\begin{cases}
	iu_{t}+\Delta u=|u|^{2k}u,\\
u(x,y,0)=f(x,y)\in H^{s}(\mathbb{R}\times \mathbb{T}).
\end{cases}
\end{equation}
and prove
\begin{thm}\label{thm: nonlineark}
Equation  \eqref{eq: knls} is {globally} well-posed for initial data in $H^{s}$, $s>1-\frac{2}{5k}$.	
\end{thm}

It is also natural to study the local well-posedness of \eqref{eq: knls} for $H^{s}$ initial data, $s=s_{c}=1-\frac{1}{k}$. Indeed, it has been proved in \cite{Wang} in the pure periodic case, with natural extension to \eqref{eq: knls}, all \eqref{eq: knls} is locally well-posedness for such initial data. The proof of \cite{Wang} relies on delicate atom space techniques, \cite{HHK,HTT}.

We want to point out, in waveguide $\mathbb{R}\times \mathbb{T}$, one can follow more standard Strichartz space time estimates scheme to obtain a relative simple LWP (i.e. local well-posedness). The $k=1$ case has already been treated in \cite{TT}. We sketch the proof for general $k$ for the convenience of readers in the Appendix B.

\subsection{Background and motivations}
\subsubsection{Bilinear estimates for Schr\"odinger}
Strichartz estimates for linear Schr\"odinger equations play fundamental role in the study of NLS. On $\mathbb{R}^{d}$, it reads, \cite{cazenave2003semilinear,keel1998endpoint},
\begin{equation}\label{eq:stri}
\|e^{it\Delta} f\|_{L_{t,x}^{\frac{2d+4}{d}}(\mathbb{R}^{d}\times \mathbb{R})}\lesssim \|f\|_{L_{x}^{2}(\mathbb{R}^{d})},
\end{equation}
while on $\mathbb{T}^{d}$, it only holds locally in time and has a $N^\epsilon$ loss, \cite{bourgain1993fourier,bourgain2014proof},
\begin{equation}
\|e^{it\Delta} P_{N}f\|_{L_{t,x}^{\frac{2d+4}{d}}(\mathbb{T}^{d}\times [0,1])}\lesssim_\epsilon N^\epsilon \|f\|_{L_{x}^{2}(\mathbb{T}^{d})}.	
\end{equation}

Let us focus on the dimension $d=2$, where $\frac{2d+4}{d}=4$, and its associated bilinear version is a $L^{2}$ estimate, which is the subject of current material.

In \cite{bourgain1998refinements}, Bourgain proves on $\mathbb{R}^{2}$,
\begin{equation}\label{eq: originalblinear}
\|e^{it\Delta} P_{N_{1}}fe^{it\Delta}P_{N_{2}}g\|_{L_{t,x}^{2}(\mathbb{R}^{2}\times \mathbb{R})}\lesssim \left(\frac{N_{2}}{N_{1}}\right)^{1/2}\|f\|_{L_{x}^{2}}\|g\|_{L_{x}^{2}}, N_{1}\geq N_{2}.	
\end{equation}
One can naturally extends his methods to generalize \eqref{eq: originalblinear} to higher dimensions, see for example \cite{colliander2008global}.

Estimate \eqref{eq: originalblinear} is useful  in the study of nonlinear problems, and in particular can be combined with high-low frequency cutoff method or I-method to establish low regularity GWP (global well-posedness), \cite{bourgain1998refinements,I-team2}.

A direct generalization of \eqref{eq: originalblinear} on $\mathbb{T}^{2}$ is much weaker than the Euclidean case, one only has, \cite{bourgain1993fourier,DPST},
\begin{equation}\label{eq: t2bilinear}
\|e^{it\Delta} P_{N_{1}}fe^{it\Delta}P_{N_{2}}g\|_{L_{t,x}^{2}(\mathbb{T}^{2}\times [0,1])}\lesssim N_{2}^{\epsilon}\|f\|_{L_{x}^{2}}\|g\|_{L_{x}^{2}}, N_{1}\geq N_{2}\geq 1.
\end{equation}
It still beats a direct application of H\"older with \eqref{eq: t2bilinear}, and plays a key role in the local well-posednessness of cubic NLS in $H^{s}(\mathbb{T}^{2})$ with $s>0$.

Motivated by implementation of I-method, bilinear estimates has also been studied in rescaled tori\footnote{We mention here, the study of bilinear estimates in rescaled tori is related to the implementation of I-method on \textbf{non-rescaled} tori.}, \cite{DPST}, it reads as
\begin{equation}\label{eq: bilinearrescaled}
\|U_{\lambda}(t) P_{N_{1}}f U_{\lambda}(t)P_{N_{2}}g\|_{L_{t,x}^{2}(\mathbb{T}_{\lambda}^{2}\times [0,1])}\lesssim N_{2}^{\epsilon}\left(\frac{1}{\lambda}+\frac{N_{2}}{N_{1}}\right)^{\frac{1}{2}}\|f\|_{L_{x}^{2}}\|g\|_{L_{x}^{2}}, N_{1}\geq N_{2}\geq 1.
\end{equation}
Here we use $U_{\lambda}(t)$ to denote the propagator of linear Schr\"odinger on $\mathbb{T}_{\lambda}^{2}$.

The proof in \cite{DPST} depends on counting, and cannot be generalized to irrational tori. Via proving a bilinear decoupling theorem following \cite{bourgain2014proof}, in \cite{FSWW}, estimates \eqref{eq: bilinearrescaled} has been extended to (rescaled) irrational tori.

In this article, we prove analogue of \eqref{eq: bilinearrescaled} on waveguide. Compared to the tori case, our main estimate \eqref{eq: bilinear1} does not have the $N_{2}^{\epsilon}$ loss, this is consistent with heruistics that waveguide behaves better than pure tori (but worse than the Euclidean case.) See also \cite{ZZ} for local-in-time bilinear estimates on waveguides.

Another point of current article is that when one has at least two directions in $\mathbb{R}$, the associated bilinear estimates \eqref{eq: bilinear2} is global in time rather than only local in time. It is in some sense no surprising. Indeed, consider the case $\mathbb{R}^{2}\times \mathbb{T}_{\lambda}^{m}$, if one does not care about the gain  $( \frac{N_2^{d-3}}{\lambda} +\frac{N_2^{d-1}}{N_1}  )^{\frac{1}{2}}$, by applying the Minkowski in all the frequency in tori directions, and Strichartz estimates in $\mathbb{R}^{2}$, one will be able to obtain some global in time estimates, (losing are large constant depending on $N_{1},N_{2},\lambda$). On the other hand, one can construct $f,g$, so that for linear Schr\"odinger equations on $\mathbb{R}\times \mathbb{T}_{\lambda}$, one has
\begin{equation}
\|U_{\lambda}(t)P_{N_{1}}fU_{\lambda}(t)P_{N_{2}}g\|_{L_{x,y,t}^{2}(\mathbb{R}\times \mathbb{T}_\lambda \times [0,\infty])}=\infty, N_{1}\gg N_{2}.	
\end{equation}
 See Appendix \ref{app: sharp} for more discussions.

 We mention that unlike the pure torus case, the analysis in the current article does not rely on the rationality of the periodic part. It may play a role if one further investigate the $\mathbb{R}\times \mathbb{T}^{n}, n\geq 2$, cases.

\subsubsection{Analysis on Waveguide}
As aforementioned, the analysis on waveguide lies somewhere between the case of torus and Euclidean cases. Our work has been in particular motivated by \cite{TT}, which proves that for linear Schr\"odinger equations on $\mathbb{R}\times \mathbb{T}$, there holds
\begin{equation}\label{eq: Stritt}
\|e^{it\Delta_{x,y}}f\|_{L_{x,y,t}^{4}(\mathbb{R}\times \mathbb{T}\times [0,1])}\lesssim \|f\|_{L_{x,y}^{2}}.
\end{equation}

The main point here is unlike the $\mathbb{T}^{2}$ case, one does not suffer from an $\epsilon$ loss in the (critical) Strichartz estimates, and thus scale invariant. See also \cite{bourgain1993fourier,killip2016scale} for scale invariant Strichartz estimates in the super-critical case.

Roughly speaking, the $L^{4}$ Strichartz estimates on (rectangular) tori, if one follows the method of counting, \cite{bourgain1993fourier}, it reduces to counting the lattice point in the region $\{\xi\in \mathbb{Z}^{2}\;|\;|\xi|^{2}= N^{2}+O(1)\}$, which ultimately reduces to counting the lattice points of a circle of radius $N$, whose numbers are bounded by $N^{\epsilon}$. The pure Euclidean case is corresponding to computing the area of $\{\xi\in \mathbb{R}^{2}\;|\;|\xi|^{2}=N^{2}+O(1)\}$, which is simple one. The waveguide case is in between, and partially makes the associated counting more close to the computation of area. We will discuss more heuristic in the proof of \eqref{eq: bilinear1}, \eqref{eq: bilinear2}.

From this viewpoint, it remains an interesting and challenging question to study $L_{x,y,t}^{\frac{2d+4}{d}}$ Strichartz estimates for waveguide with full dimension $d\geq 3.$

Barron \cite{Barron} proves global-in-time Strichartz estimates for Schr\"odinger equations on product spaces $\mathbb{R}^n \times \mathbb{T}^m$ with an additional $\epsilon$-derivative loss via $l^2$-decoupling and Euclidean Strichartz estimates. The $\epsilon$-loss is removed for exponent away from the endpoint. The estimates are applied to prove small-data-scattering at the scaling critical regularity for some nonlinear models. (See also \cite{HP} for a specific waveguide model: $\mathbb{R} \times \mathbb{T}^2$.) When $m=n=1$, the waveguide manifold is more special: in \cite{barron2020global}, the authors prove global-in-time Strichartz estimates for Schr\"odinger equations on $\mathbb{R} \times \mathbb{T}$ without a derivative loss at the endpoint. Compared to \eqref{eq: Stritt}, the main estimate in \cite{barron2020global} is stronger since it is global.

Finally, we mention there are stronger version of bilinear estimates in Euclidean spaces, which is sometimes called bilinear restriction estimates. Let us focus on the dimension $d=2$ for concreteness. Let $Q_{1}, Q_{2}$ be unit balls (or cubes) of size one in $\mathbb{R}^{2}$, which are separated by distance $\sim 1$.

One is interested to obtain estimates
\begin{equation}\label{eq: bires1}
\|e^{it\Delta}P_{Q_{1}}fe^{it\Delta}P_{Q_{2}}g\|_{L_{t,x}^{2}}\lesssim \|\hat{f}\|_{L^{p}}\|\hat{g}\|_{L^{p}}, p<2,
\end{equation}
and
\begin{equation}\label{eq: bires2}
\|e^{it\Delta}P_{Q_{1}}fe^{it\Delta}P_{Q_{2}}g\|_{L_{t,x}^{q}}\lesssim \|\hat{f}\|_{L^{2}}\|\hat{g}\|_{L^{2}}, q<2.
\end{equation}

Note that one main point of \eqref{eq: bires1}, \eqref{eq: bires2} is that $p,q<2$. And it can be applied to study concentration compactness theory, \cite{begout2007mass}.

One mainly refer to \cite{tao2003sharp,tao1998bilinear} for sharp versions of \eqref{eq: bires1}, \eqref{eq: bires2} in all dimensions.\\

It is not clear how to generalize those estimates in the waveguide case, see \cite{deng2023bilinear} for partial progress in this direction.

\subsubsection{Nonlinear problems on waveguide}
We now give a brief overview for the research line: `Nonlinear Schr\"odinger equations on waveguides'. This topic has been studied a lot in last decades since it is a hot topic in the area of nonlinear dispersive equations. The classical dispersive methods and new analysis tools are combined to investigate this topic.

The waveguide manifolds $\R^m \times \T^{d-m}$ are a product of the Euclidean space with tori, and are of particular interest in nonlinear optics of telecommunications. In fact, in today's backbone networks, data signals are almost exclusively transmitted by optical carriers in fibers (a special case of a waveguide).  Applications like the internet demand an increase in the available bandwidth in the network and a reduction of costs for the transmission of data. The nonlinear Schr\"odinger type of model is of particular importance in the description of nonlinear effects in optical fibers.

As learned from physics, an optical waveguide is a structure that `guides' a light wave by constraining it to travel along a certain desired path. One interesting feature of studying the behavior of solutions on the waveguide manifold is that it mixed inherits properties from those on classical Euclidean spaces and tori, which captures well the physics behind it. Due to the nature of such product spaces, we see NLS posed on the waveguide manifold mixed inheriting properties from those on classical Euclidean spaces and tori. 
The Euclidean case is studied and the theory, at least in the defocusing setting, is well established. (See \cite{CW,I-team1,Dodson3,Taobook} and the references therein.) Moreover, we refer to \cite{HTT1, IPT3, KV16, Yue} for a few works on tori. Due to the nature of such product spaces, we see NLS posed on the waveguide manifold mixed inheriting properties from those on classical Euclidean spaces and tori. The techniques used in Euclidean and tori settings are frequently combined and applied to waveguide problems. We refer to \cite{R2T,CGZ,CZZ,HP,HTT1,TV1,TV2,HTT2,IPT3,IPRT3,yang2023scattering,Z1,Z2,ZZ}  for some NLS results in the waveguide setting. We note that, though scattering behavior is not expected for the periodic case because of the lack of dispersive, for some specific models of waveguides, scattering results can be obtained as in the Euclidean. (See \cite{R2T,HP,TV2} for example.) 

 \cite{MR3406826} studies the asymptotic behavior of solutions to the cubic defocusing NLS posed on the waveguide manifold $\mathbb{R}\times \mathbb{T}^n$. In particular, they show that the asymptotic dynamic of small solutions is related to that of solutions of the associate resonant system and, as a consequence, they obtain global strong solutions with infinitely growing high Sobolev norms $H^s$.
\subsection{Structure of the article}
In Section 2, we present preliminary in Fourier transform, (mainly fix a notation convention), and $X^{s,b}$ space analysis. In Section 3, we present the geometric heuristic for the proof of Theorem \ref{thm: bilinearmain} and then give a rigorous analytical proof.  In Section 5, we present an improved I-method computation, proving Theorem \ref{thm: nonlinearmain} and Theorem \ref{thm: nonlineark}. In Appendix, we show the sharpness of Theorem \ref{thm: bilinearmain} and give a sketch of LWP with standard Strichartz space time estimates for \eqref{eq: knls}.

\subsection{Notations}
For waveguide, we will usually use $n$ to denote the dimension of $\mathbb{R}$, and $m$ to denote the dimension of $\mathbb{T}$ (or $\mathbb{T}_{\lambda}$), and will generally use $d=n+m$ to denote the full dimensions. We will typically use $z=(x,y)$ to denote a point in the waveguide, where $x\in \mathbb{R}^{n}$ and $y\in \mathbb{T}^{m}$.

We say $ A \lesssim B$ if there exits some $C$ so that $A\leq CB$. We say $A\gtrsim B$ if $B\lesssim A$. We say $A\sim B$ if $A \lesssim B$ and $B\lesssim A$.
We use usual Lebesgue spaces $L^{p}$ and $L_{t}^{p}L_{x}^{q}$, and Sobolev spaces $H^{s}$. In addition, $\langle a\rangle:=1+|a|$ and $a \pm:=a \pm \epsilon$ with $0<\epsilon \ll 1$.

We use $P_{\leq N}$ to denote Littlewood-Paley projection to frequency $\leq N$, and we define $P_{N}:=P_{\leq N}-P_{\leq \frac{N}{2}}$, i.e. Littlewood-Paley projection to frequency $\sim N$. A more detailed explanation is as follows.

We define the Fourier transform on $\mathbb{R}^n \times \mathbb{T}^m$ as follows:
\begin{equation}
    (\mathcal{F} f)(\xi)= \int_{\mathbb{R}^n \times \mathbb{T}^m}f(z)e^{-iz\cdot \xi}dz,
\end{equation}
where $\xi=(\xi_1,\xi_2,...,\xi_{d})\in \mathbb{R}^n \times \mathbb{Z}^m$ and $d=m+n$. We also note the Fourier inversion formula
\begin{equation}
    f(z)=c \sum_{(\xi_{n+1},...,\xi_{d})\in \mathbb{Z}^m} \int_{(\xi_1,...,\xi_{n}) \in \mathbb{R}^n} (\mathcal{F} f)(\xi)e^{iz\cdot \xi}d\xi_1...d\xi_n.
\end{equation}
For convenience, we may consider the discrete sum to be the integral with discrete measure so we can combine the above integrals together and treat them to be one integral. Moreover, we define the Schr{\"o}dinger propagator $e^{it\Delta}$ by
\begin{equation}
    \left(\mathcal{F} e^{it\Delta}f\right)(\xi)=e^{-it|\xi|^2}(\mathcal{F} f)(\xi).
\end{equation}
We are now ready to define the Littlewood-Paley projections. First, we fix $\eta_1: \mathbb{R} \rightarrow [0,1]$, a smooth even function satisfying
\begin{equation}
    \eta_1(\xi) =
\begin{cases}
1, \ |\xi|\le 1,\\
0, \ |\xi|\ge 2,
\end{cases}
\end{equation}
and $N=2^j$ is a dyadic integer. Let $\eta^d=\mathbb{R}^d\rightarrow [0,1]$, $\eta^d(\xi)=\eta_1(\xi_1)\eta_1(\xi_2)\eta_1(\xi_3)...\eta_1(\xi_d)$. We define the Littlewood-Paley projectors $P_{\leq N}$ and $P_{ N}$ by
\begin{equation}
    \mathcal{F} (P_{\leq N} f)(\xi):=\eta^d\left(\frac{\xi}{N}\right) \mathcal{F} (f)(\xi), \quad \xi \in \mathbb{R}^n \times \mathbb{Z}^m,
\end{equation}
and
\begin{equation}
P_Nf=P_{\leq N}f-P_{\leq \frac{N}{2}}f.
\end{equation}
For any $a\in (0,\infty)$, we define
\begin{equation}
    P_{\leq a}:=\sum_{N\leq a}P_N,\quad P_{> a}:=\sum_{N>a}P_N.
\end{equation}
%We use $\hat{f}$ or $\mathcal{F}f$ to denote the Fourier transform of $f$.

%We define integration on the dual of waveguide, $\mathbb{R}^{m}\times \mathbb{Z}_{1/\lambda}^{n}$ as
%\begin{equation}\label{eq: intergration}
%	\int f(\xi)(d \xi)_{\lambda}:=\frac{1}{\lambda^n} %\sum_{(\xi_{m+1},\cdots,\xi_{m+n}) %\in\mathbb{Z}_{1/\lambda}^n} \int_{(\xi_{1},\cdots,\xi_{m}) %\in \mathbb{R}^m}f(\xi) d\xi_{1}\cdots d\xi_{m}.
%\end{equation}
\section*{Acknowledgments}
The authors are grateful to the anonymous referees for their valuable comments and suggestions.  We appreciate Dr. Ze Li for helpful discussions and  beneficial suggestions on this paper.

C. Fan was supported by the National Key R\&D Program of China, 2021YFA1000800, CAS Project for Young Scientists in Basic Research, Grant No. YSBR-031, and NSFC
Grant No. 12288201.

 K. Yang was supported by Chongqing Science and Technology Commission (ncamc2022-msxm04), the NSF grant of China (No.12301296, 12371244) and the Science and Technology Research Program of Chongqing Municipal Education Commission
(Grant No.KJQN202200506).

 Z. Zhao was supported by the NSF grant of China (No. 12101046, 12271032) and the Beijing Institute of Technology Research Fund Program for Young Scholars.
 
  J. Zheng was supported by National key R\&D program of China: 2021YFA1002500, NSF grant of China (No. 12271051) and Beijing Natural Science Foundation 1222019.

\section{Preliminary}
\subsection{Fourier Transform in the waveguide setting}
We quickly recall the Fourier transform in the (rescaled) waveguide case. There seems to be different conventions of (equivalent) definitions. We mainly follow the presentation in \cite{DPST}.\\

For notation convenience, we will sometimes treat $\mathbb{T}$ as $[0,1]$, (and accordingly, $\mathbb{T}_{\lambda}$ as $[0,\lambda]$).

Let Fourier transform on $\mathbb{R}^m \times \mathbb{T}_\lambda^n$ be defined as follows,
\begin{equation}
 \hat{f}(\xi)= (\mathcal{F} f)(\xi)= \int_{\mathbb{R}^m \times \mathbb{T}_\lambda^n}f(z)e^{-2\pi i z\cdot \xi}dz,
\end{equation}
where $z=(x,y)\in \mathbb{R}^m \times \mathbb{T}_\lambda^n$ and $\xi=(\xi_1,\xi_2,...,\xi_{m+n})\in \mathbb{R}^m \times \mathbb{Z}_{1/\lambda}^n$.

If one further defines the integration\footnote{Part of the integration is really discrete summations.} on $\mathbb{R}^m \times \mathbb{Z}_{1/\lambda}^n$ as
\begin{equation}
	\int f(\xi)(d \xi)_{\lambda}:=\frac{1}{\lambda^n} \sum_{(\xi_{m+1},\cdots,\xi_{m+n}) \in\mathbb{Z}_{1/\lambda}^n} \int_{(\xi_{1},\cdots,\xi_{m}) \in \mathbb{R}^m}f(\xi) d\xi_{1}\cdots d\xi_{m}.
\end{equation}

With those notations, one has, as usual,

\begin{enumerate}
	\item Fourier inversion formula
	 \begin{equation}
    f(z)= \frac{1}{\lambda^n} \sum_{(\xi_{m+1},...,\xi_{d})\in \mathbb{Z}_{1/\lambda}^n} \int_{(\xi_1,...,\xi_{m}) \in \mathbb{R}^m} (\mathcal{F} f)(\xi)e^{2\pi iz\cdot \xi}d\xi_1...d\xi_m.
\end{equation}
    \item Plancherel identity
    \begin{equation}
    \|f\|_{L^2\left(\mathbb{R}^m \times \mathbb{T}_\lambda^n\right)}=\|\hat{f}\|_{L^2\left(\mathbb{R}^m \times \mathbb{Z}_{1/\lambda}^n\right)  }.  \end{equation}
    \item Convolution and multiplication under Fourier transform
    \begin{equation}
    	\widehat{f g}(\xi)=\hat{f} \star \hat{g}(\xi):=\int_{\mathbb{R}^m \times \mathbb{Z}_{1/\lambda}^n} \hat{f}\left(\xi-\eta\right) \hat{g}\left(\eta\right)\left(d \eta\right)_\lambda.
    \end{equation}

\end{enumerate}

\subsection{$X^{s,b}$ spaces and transference principle}
Following \cite{bourgain1993fourier,bourgain1998refinements}, one may define $X^{s,b}$ norm as
\begin{equation}
	\|u\|_{X^{s,b}}:=\| \langle \tau-|\xi|^{2}\rangle^{b}\langle \xi\rangle^{s}\tilde{u}\|_{L_{\xi,\tau}^{2}},
\end{equation}
where $u=u(z,t)$ is a function defined on $\mathbb{R}^{m}\times \mathbb{T}_{\lambda}^{n}\times \mathbb{R}$, and $z\in \mathbb{R}^{m}\times \mathbb{T}_{\lambda}^{n}, t\in \mathbb{R}$. And $\tilde{u}(\xi,\tau)$ is the space-time Fourier transform of $u$, where $\xi\in \mathbb{R}^{m}\times \mathbb{Z}_{1/\lambda}^{n}$, $\tau\in \mathbb{R}$.

And $X^{s,b}$ spaces are just all those functions with finite $X^{s,b}$ norm.

In practice, one mainly works on $s\geq 0,b>\frac{1}{2}$.

 One key property for $X^{s,b}$ space is that inherits the estimates of linear solutions, \cite{bourgain1998refinements}.  In particular, Theorem \ref{thm: bilinearmain} implies
\begin{lem}\label{lem: trans}
Let $b>\frac{1}{2}$, and $u_{1}, u_{2}\in X^{0,b}(\mathbb{R}^{m}\times \mathbb{T}_{\lambda}^{n}\times (-\infty,\infty))$.

For $n=1,m=1$, one has
\begin{equation}\label{eq: xsbbilinear1}
 \|  P_{N_1}u_{1}P_{N_2}u_{2} \|_{L_{x,y,t}^2(\mathbb{R}\times \mathbb{T}_{\lambda} \times [0,1])} \lesssim \left(\frac{1}{\lambda} +\frac{N_2}{N_1}  \right)^{\frac{1}{2}} \|u_{1}\|_{X^{0,b}} \|u_{2}\|_{X^{0,b}}.
\end{equation}

For $n\geq 2, m\geq 1$, one has
\begin{equation}\label{eq: xsbbilinear2}
 \| P_{N_1}u_{1}P_{N_2}u_{2}\|_{L_{x,y,t}^2(\mathbb{R}^{m}\times \mathbb{T}^{n}_{\lambda} \times [0,\infty))} \lesssim \left( \frac{N_2^{d-3}}{\lambda} +\frac{N_2^{d-1}}{N_1}  \right)^{\frac{1}{2}}  \|u_{1}\|_{X^{0,b}} \|u_{2}\|_{X^{0,b}},
\end{equation}
where $N_{1}\geq N_{2}\geq 1$.
\end{lem}

\section{Proof of Theorem \ref{thm: bilinearmain}} \label{the proof of bilinear}
\subsection{Geometric Heuristic}
In this subsection, we first briefly explain the geometric heuristics for Strichartz estimates on $\mathbb{T}^{2}$\cite{bourgain1993fourier}, and  $\mathbb{R}\times \mathbb{T}$ \cite{TT}, and then illustrate the heuristics for our estimates \eqref{eq: bilinear1}.

Let $C_{1,N}, C_{2,N}$ be the optimal constant for $L_{t,z}^{4}$ Strichartz estimates on  $\mathbb{T}^{2}$, and  $\mathbb{R}\times \mathbb{T}$, respectively.  Thus
\begin{align}\label{eq: s1}
\|e^{it\Delta}P_{N}f\|_{L_{t,z}^{4}(\mathbb{T}^{2}\times [0,1])}\leq& C_{1,N}\|f\|_{L_{z}^{2}}\\\label{eq: s2}
	\|e^{it\Delta}P_{N}f\|_{L_{t,z}^{4}(\mathbb{R}\times \mathbb{T}\times [0,1])}\leq& C_{2,N}\|f\|_{L_{z}^{2}}.
\end{align}

The bound of $C_{1,N}$ is essentially counting the numbers $(n_{1}, n_{2})$ in $\mathbb{Z}^{2}$, such that
\begin{equation}
n_{1}^{2}+n_{2}^{2}=N^{2}+O(1)	
\end{equation}

Note that this gives an annulus of radius $N$, but thickness $\sim \frac{1}{N}$. Thus, if the counting of numbers is like computing area, this should gives a bound $\sim 1$.
But one cannot exclude the possibility that there are many lattice points on some circle\footnote{This also explains why a direct counting argument on irrational tori cannot work, since one needs to count, for example, $n_{1}^{2}+\gamma n_{2}^{2}=N^{2}+O(1)$, $\gamma$ irrational. In this case, one cannot reduce such a counting to finite many circles $n_{1}^{2}+\gamma n_{2}^{2}=R, R=N^{2}+O(1)$.} $n_{1}^{2}+n_{2}^{2}=R$, $R=N^{2}+O(1)$.

In the situation of waveguide, to bound $C_{2,N}$, one is essentially counting the total length of all the intervals which is the intersection of the line $y=n$, for some $n\in \mathbb{Z}$, with the region $x^{2}+y^{2}=N^{2}+O(1)$, see Figure 1 below.
\begin{center}
\begin{tikzpicture}
\tkzDefPoints{0/0/O,2.5/0/A,0/{1/3}/O_1,0/{2/3}/O_2, 0/1/O_3, 0/{4/3}/O_4, 0/{5/3}/O_5, 0/2/O_6, 0/{-1/3}/O_7, 0/{-2/3}/O_8, 0/-1/O_9, 0/{-4/3}/O_{10}, 0/{-5/3}/O_{11}, 0/-2/O_{12}, 0/{-7/3}/O_{13}, 0/{7/3}/O_{14}}
\tkzDefPointWith[linear, K=9/10](O,A) \tkzGetPoint{A_1}
\tkzDefLine[parallel=through O_1](O,A) \tkzGetPoint{l_1}
\tkzDefLine[parallel=through O_2](O,A) \tkzGetPoint{l_2}
\tkzDefLine[parallel=through O_3](O,A) \tkzGetPoint{l_3}
\tkzDefLine[parallel=through O_4](O,A) \tkzGetPoint{l_4}
\tkzDefLine[parallel=through O_5](O,A) \tkzGetPoint{l_5}
\tkzDefLine[parallel=through O_6](O,A) \tkzGetPoint{l_6}
\tkzDefLine[parallel=through O_7](O,A) \tkzGetPoint{l_7}
\tkzDefLine[parallel=through O_8](O,A) \tkzGetPoint{l_8}
\tkzDefLine[parallel=through O_9](O,A) \tkzGetPoint{l_9}
\tkzDefLine[parallel=through O_{10}](O,A) \tkzGetPoint{l_{10}}
\tkzDefLine[parallel=through O_{11}](O,A) \tkzGetPoint{l_{11}}
\tkzDefLine[parallel=through O_{12}](O,A) \tkzGetPoint{l_{12}}
\tkzDefLine[parallel=through O_{13}](O,A) \tkzGetPoint{l_{13}}
\tkzDefLine[parallel=through O_{14}](O,A) \tkzGetPoint{l_{14}}
\begin{scope}
\tkzClipCircle(O,A)
\tkzFillCircle[color=gray!20](O,A)
\tkzDrawLine[add= 5 and 1](O,A_1)
\tkzDrawLine[add= 1 and 1](O_1,l_1)
\tkzDrawLine[add= 1 and 1](O_2,l_2)
\tkzDrawLine[add= 1 and 1](O_3,l_3)
\tkzDrawLine[add= 1 and 1](O_4,l_4)
\tkzDrawLine[add= 1 and 1](O_5,l_5)
\tkzDrawLine[add= 1 and 1](O_6,l_6)
\tkzDrawLine[add= 1 and 1](O_7,l_7)
\tkzDrawLine[add= 1 and 1](O_8,l_8)
\tkzDrawLine[add= 1 and 1](O_9,l_9)
\tkzDrawLine[add= 1 and 1](O_{10},l_{10})
\tkzDrawLine[add= 1 and 1](O_{11},l_{11})
\tkzDrawLine[add= 1 and 1](O_{12},l_{12})
\tkzDrawLine[add= 1 and 1](O_{13},l_{13})
\tkzDrawLine[add= 1 and 1](O_{14},l_{14})
\tkzFillCircle[color=white](O,A_1)
\end{scope}
\tkzDrawCircles(O,A_1  O,A)
\tkzLabelPoints[below](O)
\tkzDrawPoints(O)
\node at (0,-3){Figure 1};
\end{tikzpicture}
\end{center}

Note that the vertical distance of those intervals are 1, thus the total length of those intervals are essentially the area of the annulus, (which is $\sim 1$), with an error bounded by the length of longest inverval.

But it is elementary geometry for an annulus of radius $N$, and thickness $\frac{1}{N}$, the longest length of such interval can be at most $\sqrt{N\times \frac{1}{N}}\sim 1$.

This has been made rigorous in \cite{TT} and bounds $C_{2,N}$ by a constant independent of $N$.

Now, let $C_{3,N_{1},N_{2}}$ the optimal constant such that
\begin{equation}\label{eq: b1}	\| U_{\lambda}(t) P_{N_1}f  U_{\lambda}(t) P_{N_2}g \|_{L_{x,y,t}^2(\mathbb{R}\times \mathbb{T}_{\lambda} \times [0,1])} \lesssim C_{3,N_{1},N_{2}}^{\frac{1}{2}} \|f\|_{L^2} \|g\|_{L^2}
\end{equation}

Then, one may follow \cite{bourgain1993fourier,bourgain1998refinements}, to reduce the bound of $C_{3,N_{1},N_{2}}$ to  $\frac{1}{\lambda}$ total length of all the intervals in the shaded area in Figure 2 below.
\begin{center}
\begin{tikzpicture}
\tkzDefPoints{0/0/O,-{9/10}/{-12/10}/a,-6/-6/c, 1/0/e, 0/{1/5}/i, -3/-3/A, -3/3/B, 3/3/C, 3/-3/D}
\tkzDefPointBy[translation= from O to e](a)  \tkzGetPoint{d}
\tkzDefPointBy[translation= from O to i](d)  \tkzGetPoint{d_1}
\tkzDefLine[parallel=through d_1](a,d) \tkzGetPoint{l_1}
\tkzDefPointBy[translation= from i to O](d)  \tkzGetPoint{d_{-1}}
\tkzDefLine[parallel=through d_{-1}](a,d) \tkzGetPoint{l_{-1}}
\tkzDefPointBy[translation= from O to i](d_1)  \tkzGetPoint{d_2}
\tkzDefLine[parallel=through d_2](a,d) \tkzGetPoint{l_2}
\tkzDefPointBy[translation= from O to i](d_2)  \tkzGetPoint{d_3}
\tkzDefLine[parallel=through d_3](a,d) \tkzGetPoint{l_3}
\tkzDefPointBy[translation= from O to i](d_3)  \tkzGetPoint{d_4}
\tkzDefLine[parallel=through d_4](a,d) \tkzGetPoint{l_4}
\tkzDefPointBy[translation= from O to i](d_4)  \tkzGetPoint{d_5}
\tkzDefLine[parallel=through d_5](a,d) \tkzGetPoint{l_5}
%\tkzDefMidPoint(a,b) \tkzGetPoint{c}
\tkzDefPointWith[linear, K=53/50](c,a) \tkzGetPoint{a_1}
\tkzDefPointWith[linear, K=51/50](c,a) \tkzGetPoint{a_2}
\begin{scope}

\tkzClipCircle(O,a)
\tkzClipCircle(c,a_1)
\tkzFillCircle[color=gray!20](c,a_1)
\tkzDrawLine[add= 10 and 10](a,d)
\tkzDrawLine[add= 10 and 10](d_1,l_1)
\tkzDrawLine[add= 10 and 10](d_2,l_2)
\tkzDrawLine[add= 10 and 10](d_3,l_3)
\tkzDrawLine[add= 10 and 10](d_4,l_4)
\tkzDrawLine[add= 10 and 10](d_5,l_5)
\tkzDrawLine[add= 10 and 10](d_{-1},l_{-1})

\tkzFillCircle[color=white](c,a_2)
\end{scope}

\begin{scope}
\tkzClipPolygon(A,B,C,D)
\tkzDrawArc[rotate,color=red](c,a_1)(30)
\tkzDrawArc[rotate,color=red](c,a_1)(-30)
\tkzDrawArc[rotate,color=red](c,a_2)(30)
\tkzDrawArc[rotate,color=red](c,a_2)(-30) 
\end{scope}

\tkzDrawCircles[color=blue](O,a)
\tkzLabelPoints[below](O)
\tkzDrawPoints(O)
\node at (0,-4){Figure 2};
\end{tikzpicture}

\end{center}

Here, $O$ is original point, the blue part is a circle of radius $N_{2}$, centered at origin. The red part is an annulus of radius $N_{1}$, thickness $\frac{1}{N_{1}}$. And the intervals in the shaded area are obtained by intersect line $y=\frac{n}{\lambda}$, for some $n\in Z$ with the annulus, (we use $(x,y)$ to denote a point in the plane.) Note that the vertical distance of those intervals are $\frac{1}{\lambda}$.

The bound of $C_{3,N_{1},N_{2}}$ basically follows from three facts,
\begin{enumerate}
	\item The area of shaded part is bound by $\frac{N_{2}}{N_{1}}$.
	\item The total length of those intervals multiplying $\frac{1}{\lambda}$ is bounded by the area of shaded part with an error bounded by $\frac{1}{\lambda}$ multiplying longest length of such intervals.
	\item As aforementioned, for an annulus of radius $N_1$, and thickness $\frac{1}{N_1}$, the longest length of such interval can be at most $\sqrt{N_1\times \frac{1}{N_1}}\sim 1$.
\end{enumerate}

This will give our desired bound.

For the case $\bR^2 \times \bT_\lambda$ restricted to time interval $[0,1]$,  we need to bound the $\frac{1}{\lambda}$ total areas of the all intersections of  planes which are of the type $\{z=\frac{k}{\lambda}\} $, $k\in \bZ$, a ball which is of radius $N_2$ and a spherical shell which is of radius $N_1$ and thickness $\frac{1}{N_1}$. Similarly, the largest area of such intersections is at most $1$, and the volume of the intersection of the ball and the spherical shell is $\sim \frac{N_2^2}{N_1}$.

When we consider the case $m\ge 2, n\ge 1$, we will fix $d-3$ components, and reduce the counting problem to $\bR^2\times \bT_\lambda$ case. Each component contributes $N_2$ choices, so the final bound is $(\frac{1}{\lambda}+\frac{N_2^2}{N_1})N_2^{d-3}$. If we restrict to the long time interval $[0,T]$, the analysis is similar, and we will see that why we need at least two directions in $\bR$.

We mention that the parallel counting problems in the pure torus case has been studied in \cite{DPST}.

We now turn to more rigorous analysis.
\subsection{Proof of Estimate \eqref{eq: bilinear1}}
When $N_{1}\sim N_{2}$, estimate \eqref{eq: bilinear1} follows from  Strichartz estimate on $\bR \times \bT$ in
\cite{TT}.

We now focus on the case $N_1\gg N_2$. 

We follow the argument in the proof of Proposition 3.7 in \cite{DPST} to reduce \eqref{eq: bilinear1} to a measure estimate problem (or counting problem).

Let $\phi(\tau)$ be some nice function with $\hat{\phi}(t)\gtrsim 1$ on [-1,1], and $\supp \phi\subset [-\frac{1}{2}, \frac{1}{2}]$.

By Plancharel's theorem,
\begin{align}
 & \| U_\lambda(t) P_{N_1}f \cdot U_\lambda(t) P_{N_2}g      \|_{L^2(\bR \times \bT_\lambda\times [0,1])}  \nonumber  \\
 \lesssim & \| \left(\widehat{\phi}(t) U_\lambda(t) P_{N_1}f \right) \cdot \left( \widehat{\phi}(t) U_\lambda(t) P_{N_2}g \right)  \|_{L^2_t(\bR) L^2_z(\bR \times \bT_{\lambda})} \nonumber \\
 = &   \left\| \int_{\eta=\eta^1+\eta^2, \atop \tau=\tau^1+\tau^2} \phi(\tau^1-|\eta^1|^2) \phi(\tau^2-|\eta^2|^2) \widehat{P_{N_1}f}(\eta^1) \widehat{P_{N_2}g}(\eta^2)     ({\rm d}\eta^1)_{\lambda} ({\rm d}\eta^2)_{\lambda} {\rm d}\tau^1 {\rm d}\tau^2   \right\|_{L^2_\tau(\bR) L^2_\eta(\bR \times \bZ_{1/\lambda})}   \nonumber \\
 = &  \left\| \int_{\eta=\eta^1+\eta^2} \widetilde{ \phi}(\tau-|\eta^1|^2-|\eta^2|^2) \widehat{P_{N_1}f}(\eta^1) \widehat{P_{N_2}g}(\eta^2)     ({\rm d}\eta^1)_{\lambda} ({\rm d}\eta^2)_{\lambda}    \right\|_{L^2_\tau(\bR) L^2_\eta(\bR \times \bZ_{1/\lambda})},  \label{th:bilinear estimate-eq-1}
\end{align}
where $\widetilde{\phi}\in \mathcal{S}$ is defined by
$$ \int_{\tau=\tau^1+\tau^2} \phi(\tau^1-|\eta^1|^2) \phi(\tau^2-|\eta^2|^2) {\rm d}\tau^1 {\rm d}\tau^2 = \widetilde{ \phi}(\tau-|\eta^1|^2-|\eta^2|^2),$$
and it is easy to check that $\supp(\widetilde{\phi})\subset [-1, 1]$.

Using Cauchy-Schwartz and the definition of $P_{N}$,
\begin{align*}
& \left|\int_{\eta=\eta^1+\eta^2} \widetilde{ \phi}(\tau-|\eta^1|^2-|\eta^2|^2) \widehat{P_{N_1}f}(\eta^1) \widehat{P_{N_2}g}(\eta^2)     ({\rm d}\eta^1)_{\lambda} ({\rm d}\eta^2)_{\lambda} \right| \\
\lesssim & M(\eta,\tau) \left(\int_{\eta=\eta^1+\eta^2} \widetilde{ \phi}(\tau-|\eta^1|^2-|\eta^2|^2) |\widehat{f}(\eta^1)|^2 |\widehat{g}(\eta^2)|^2     ({\rm d}\eta^1)_{\lambda} ({\rm d}\eta^2)_{\lambda} \right)^{1/2},
\end{align*}
where
\begin{align*}
M(\eta,\tau) & =  \left( \int_{\eta=\eta^1+\eta^2} \widetilde{ \phi}(\tau-|\eta^1|^2-|\eta^2|^2) \left|  \widehat{P_{N_1}\delta_0} \right|^2 \left|  \widehat{P_{N_2}\delta_0} \right|^2 ({\rm d}\eta^1)_{\lambda} ({\rm d}\eta^2)_{\lambda}  \right)^{1/2} \\
& \lesssim\left( \int_{\eta=\eta^1+\eta^2, \atop |\eta^1|\sim N_1, |\eta^2|\sim N_2 }{\bf 1}_{\{|\tau-|\eta^1|^2-|\eta^2|^2|\le 1\}}(\tau, \eta^1, \eta^2 ) ({\rm d}\eta^1)_{\lambda} ({\rm d}\eta^2)_{\lambda}  \right)^{1/2}.
\end{align*}
So we could bound (\ref{th:bilinear estimate-eq-1}) by
\begin{align*}
& \left\|M(\eta,\tau) \left(\int_{\eta=\eta^1+\eta^2} \widetilde{ \phi}(\tau-|\eta^1|^2-|\eta^2|^2) |\widehat{f}(\eta^1)|^2 |\widehat{g}(\eta^2)|^2     ({\rm d}\eta^1)_{\lambda} ({\rm d}\eta^2)_{\lambda} \right)^{1/2} \right\|_{L^2_\tau L^2_\eta} \\
\lesssim &\|M(\eta,\tau)\|_{L^\infty_\tau L^\infty_\eta} \left\|\left(\int_{\eta=\eta^1+\eta^2} \widetilde{ \phi}(\tau-|\eta^1|^2-|\eta^2|^2) |\widehat{f}(\eta^1)|^2 |\widehat{g}(\eta^2)|^2     ({\rm d}\eta^1)_{\lambda} ({\rm d}\eta^2)_{\lambda} \right)^{1/2} \right\|_{L^2_\tau L^2_\eta} \\
\lesssim & \|M(\eta,\tau)\|_{L^\infty_\tau L^\infty_\eta} \|f\|_{L^2} \|g\|_{L^2}.
\end{align*}
To get
$$\|M(\eta,\tau)\|_{L^\infty_\tau L^\infty_\eta}\lesssim \left( \frac{1}{\lambda}+\frac{N_2}{N_1} \right)^{1/2},$$
it suffices to prove the measure estimate
\begin{equation}\label{eq: measure estimate}
\left| \{ \xi\in \bR \times \bZ_{1/\lambda} : |\xi|\sim N_2, \left| |\xi|^2+|\eta-\xi|^2-\tau \right|\le 1 \}  \right| \lesssim \frac{1}{\lambda}+\frac{N_2}{N_1}
\end{equation}
for any fixed $|\eta|\sim N_1, \tau\in \bR$.

We now fix $\eta=(\eta_1,\eta_2) \in \bR \times \bZ_{1/\lambda}, \tau\in \bR$ with $|\eta|\sim N_1$, and define
$$ C=\{ \xi\in \bR \times \bZ_{1/\lambda} : |\xi|\sim N_2, \left| |\xi|^2+|\eta-\xi|^2-\tau \right| \le 1  \}.$$

Let's first analyze the value range for the real or discrete component of $\xi$ when the other component is fixed. We claim that when $|\eta_1|\sim N_1$, \eqref{eq: measure estimate} holds.  Fix $\xi\p=(\xi\p_1, \xi\p_2) \in C$, so if $\xi=(\xi_1, \xi\p_2) \in C$, then
$$ \left||\xi|^2+|\eta-\xi|^2-|\xi\p|^2-|\eta-\xi\p|^2   \right|\le 2,$$
that is
$$ |(\xi_1-\xi\p_1)(\xi_1+\xi\p_1-\eta_1)| \le 1 , $$
so
$$ |\xi_1-\xi\p_1| \lesssim \frac{1}{N_1}.$$
Thus
\begin{equation}
\begin{aligned}
&\int_{  \bR \times \bZ_{1/\lambda} } 1_{C}(\xi) \dif \xi \\
\lesssim & \int_{|\xi_2 |\lesssim N_2}    \int 1_{\{\xi_1: \xi \in C  \} }(\xi_1)  \dif \xi\\
\lesssim &\frac{N_2}{N_1}. 
\end{aligned}
\end{equation}

The above arguments also holds when $|\eta_2|\sim N_1$ and  $\lambda \gtrsim N_1$. Now we suppose that $\lambda \ll N_1$, and we assume $\eta_2>0$ and $\eta_2\sim N_1$  without loss of generality.

We hope to clearly indicate the value range of $\xi_1$ when $\xi_2=\frac{k}{\lambda}\in \bZ_{1/\lambda}$ is fixed. For this purpose, we define
$$ \mu_k= -\left|\frac{k}{\lambda}-\frac{\eta_2}{2}\right|^2-\frac{|\eta|^2}{4}+\frac{\tau}{2}   $$
for $k \in \bZ$, then
$$C=\l(\bigcup_{k\in \bZ} C_k\r) \bigcap \{|\xi|\sim N_2\},$$
where
$$ C_k= \l\{ \xi \in \bR \times \bZ_{1/\lambda}: \l| \l|\xi_1-\frac{\eta_1}{2} \r|^2-\mu_k \r|\le \frac12, \xi_2=\frac{k}{\lambda}\r\},$$
so we only need to consider $k$ that satisfies $|k|\lesssim \lambda N_2$ and $\mu_k > -\frac12$ for estimating $|C|$.

Through direct calculation and $\eta_2\sim N_1\gg N_2$, we see that $\{\mu_k\}_{|k|\lea \lambda N_2}$ is almost an arithmetic sequence, that is
$$ \mu_{k+1}-\mu_k \sim \frac{N_1}{\lambda}\gg 1,$$
so there are at most $O(1)$ $k$'s such that $|\mu_k|\le 1$, then
$$ \l|\bigcup_{k:|\mu_k|\le 1}C_k \r| \lesssim \frac{1}{\lambda}.$$

Now we only consider $k$'s that satisfy  $\mu_k> 1$, for these $k$'s,
$$ C_k=\l\{ \xi \in \bR \times \bZ_{1/\lambda}: (\mu_k- \frac12)^{\frac{1}{2}} \le \l| \xi_1-\frac{\eta_1}{2} \r|\le (\mu_k+ \frac12)^{\frac{1}{2}}, \xi_2=\frac{k}{\lambda}\r\},$$
thus we can directly estimate
$$ |C_k|\sim \lambda^{-1}\mu_k^{-\frac12}.$$
There won't be many $k$'s due to $|\xi_1|\lea N_2$.     Let
$$ k_0=\min\{k\in \bZ:   C_k \bigcap \{|\xi|\sim N_2\}\ne \varnothing, \mu_k >1 , |k|\lesssim \lambda N_2       \},$$
and
$$ k_1=\max\{k\in \bZ:   C_k \bigcap \{|\xi|\sim N_2\}\ne \varnothing, \mu_k >1 , |k|\lesssim \lambda N_2       \}.$$
If $k_0=k_1$,  then we have
\begin{equation}\label{eq:estimate when k0 is k1}
\sum_{k=k_0}^{k_1} |C_k| \lesssim \lambda^{-1} \mu_{k_0}^{-\frac12} \lesssim \lambda^{-1}.
\end{equation}
 Now we assume $k_1>k_0$.   Then there holds
$$ \left[(\mu_{k_i}-\frac12)^{\frac12}, (\mu_{k_i}+\frac12)^{\frac12}  \right] \bigcap \l\{ \l|\xi_1- \frac{\eta_1}{2}\r|: |\xi_1| \lesssim N_2   \r\} \ne \varnothing,  \qquad i=0,1,$$
thus
$$  \left(\mu_{k_1}-\frac12\right)^{\frac12}-\left(\mu_{k_0}+\frac12\right)^{\frac12} \lesssim N_2,$$
so
\begin{equation}\label{eq:estimate k1-k0}
\mu_{k_1} -\mu_{k_0} \lesssim N_2 \mu_{k_1}^{\frac12}.
\end{equation}

When $\mu_{k_0} \lesssim \mu_{k_1}-\mu_{k_0}$, then we get
$$ k_1-k_0 \lesssim \frac{\lambda N_2^2}{N_1},$$
so we could estimate
\begin{align}
\sum_{k=k_0}^{k_1} |C_k| & \lesssim  \sum_{k=k_0}^{k_1} \lambda^{-1}\mu_k^{-\frac12} \lesssim \sum_{k=k_0}^{k_1} \lambda^{-1}\left(\mu_{k_0}+ (k-k_0) \frac{N_1}{\lambda}  \right)^{-\frac12} \label{eq:estimate when muk0 small}  \\
&\lesssim \frac{1}{\lambda} + \left(\frac{1}{\lambda N_1}\right)^{\frac12} (k_1-k_0)^{\frac12}      \lesssim \frac{1}{\lambda}+  \frac{N_2}{N_1}. \nonumber
\end{align}
When $\mu_{k_0} \gg \mu_{k_1}-\mu_{k_0}$, by (\ref{eq:estimate k1-k0}), we get
$$ k_1-k_0 \lesssim \frac{\lambda N_2}{N_1} \mu_{k_0}^{\frac12},$$
thus
\begin{align}
  \sum_{k=k_0}^{k_1} |C_k| & \lesssim  \sum_{k=k_0}^{k_1} \lambda^{-1}\mu_k^{-\frac12} \lesssim \lambda^{-1} (k_1-k_0+1) \mu_{k_0}^{-\frac12} \label{eq:estimate when muk0 big}    \\
 & \lesssim \frac{1}{\lambda} +\frac{N_2}{N_1}. \nonumber
\end{align}
The proof of \eqref{eq: bilinear1} is now complete.

\hfill$\Box$\vspace{2ex}

\subsection{Proof of \eqref{eq: bilinear2}}
%\noindent{\em Proof of Theorem \ref{th:global bilinear estimate}.} Now
It suffices to show that, for any $T>1$, there holds
$$ \| U_\lambda(t) P_{N_1}f \cdot U_\lambda(t) P_{N_2}g   \|_{L^2(\bR^m \times \bT^n_\lambda\times [-T,T])} \lesssim \l( K(\lambda, N_1, N_2)  \r)^{\frac12} \|f\|_{L^2} \|g\|_{L^2}.$$

Similar as the proof of \eqref{eq: bilinear1}, we only need to prove the measure estimate, that is, for fixed $|\eta|\sim N_1$, $\tau\in \bR$,
$$ T\int_{  \bR^m \times \bZ^n_{1/\lambda} } 1_{C}(\xi) \dif \xi \lesssim \l( \frac{1}{\lambda T} +\frac{N_2^2}{ N_1 T}  \r)N_2^{d-3},$$
and equivalently,
$$ |C|\lesssim  \l( \frac{1}{\lambda T} +\frac{N_2^2}{ N_1 T}  \r)N_2^{d-3},$$
where
$$ C=\left\{ \xi\in \bR^m \times \bZ^n_{1/\lambda} : |\xi|\sim N_2, \left| |\xi|^2+|\eta-\xi|^2-\tau \right| \le \frac{1}{T} \right \}.$$

The later argument is very similar to the previous one in  the proof of \eqref{eq: bilinear1}, the difference is that, we will reduce the estimate to  $\bR \times \bZ_{\lambda}$ case with some more stronger assumptions, and so the proof relies on $m\ge 2$. We will clearly see this in Case 2(a).

{\bf Case 1.}  $N_1\sim N_2$. It suffices to show that
$$ |C|\lesssim \frac{N_2^{d-2}}{T}.$$
Note that $m \ge 2$. it is easy to see that
\begin{equation}\label{eq:estimate area directly}
\l|\left\{(\xi_1, \xi_2) \in \bR^2: \l| |\xi_1-\frac{\eta_1}{2}|^2+| \xi_2-\frac{\eta_2}{2}|^2-c \r| \le \frac{1}{2T}   \right\} \r| \lesssim \frac{1}{T}
\end{equation}
holds uniformly for all $c \in \bR$, thus
\begin{align*}
|C| &\lesssim \int_{ |(\xi_3,\cdots, \xi_d)|\lesssim N_2 } \int_{(\xi_1,\xi_2)\in \bR^2} 1_C(\xi) \dif \xi \lesssim \frac{N_2^{d-2}}{T}.
\end{align*}

{\bf Case 2.}  $N_1\gg N_2$.  Similarly, we only need to consider the case when $\lambda \ll N_1 T$, and we assume that $\eta_d >0$ and $\eta_d\sim N_1\gg N_2$. In this time, $\eta_1, \cdots, \eta_m \in \bR$ and $\eta_{m+1}, \cdots, \eta_d \in \bZ_{1/\lambda}$.

{\bf Case 2(a).} $m= 2, n=1$.    Define
$$ \mu_k= -\left|\frac{k}{\lambda}-\frac{\eta_3}{2}\right|^2-\frac{|\eta|^2}{4}+\frac{\tau}{2}$$
for $k \in \bZ$,  then
$$C=\l(\bigcup_{k\in \bZ} C_k\r) \bigcap \{|\xi|\sim N_2\},$$
where
$$ C_k= \l\{ \xi \in \bR^2 \times \bZ_{1/\lambda}: \l| \l|\xi_1-\frac{\eta_1}{2} \r|^2+\l|\xi_2-\frac{\eta_2}{2} \r|^2-\mu_k \r|\le \frac{1}{2T}, \xi_3=\frac{k}{\lambda}\r\}.$$
Like (\ref{eq:estimate area directly}),    $|C_k|$ has a bound that $|C_k|\lesssim \frac{1}{\lambda T}$. We only consider $|k|\lesssim \lambda N_2$ and we have
$$ \mu_{k+1}-\mu_k \sim \frac{N_1}{\lambda}\gg \frac{1}{T},$$
then
$$ \l|\bigcup_{k:|\mu_k|\le 10\frac{N_1}{\lambda}}C_k \r| \lesssim \frac{1}{\lambda T}.$$

Now we only consider $k$'s that satisfy  $\mu_k> 10\frac{N_1}{\lambda}$.  Note that at this time, the following holds
$$ C_k= C_k \bigcap \l(\bigcup_{i=1}^2 \l\{\xi\in \bR^2 \times \bZ_{1/\lambda}:  \l|\xi_i-\frac{\eta_i}{2} \r|^2 >4 \frac{N_1}{\lambda}      \r\}      \r),$$
by symmetry, we only need to consider $ C_k \bigcap \l\{\xi\in \bR^2 \times \bZ_{1/\lambda}:  \l|\xi_1-\frac{\eta_1}{2} \r|^2 >4 \frac{N_1}{\lambda}      \r\}$.

 Fix $\xi_2$ with $|\xi_2|\lesssim N_2$, now how we analyze the value range of $(\xi_1,\xi_3)$ is similar to $m=1,n=1$ case, but there is a new restriction $\l|\xi_1-\frac{\eta_1}{2} \r|^2 >4 \frac{N_1}{\lambda}$, so we will get a better estimate.

Define
$$ \mu_k\p(\xi_2)= \mu_k- \l|\xi_2-\frac{\eta_2}{2} \r|^2, $$
\begin{equation*}
    \begin{split}
          C_k\p(\xi_2)= \Big\{ (\xi_1,\xi_3)&\in \bR \times \bZ_{1/\lambda}:  \l| \l|\xi_1-\frac{\eta_1}{2} \r|^2 -\mu_k\p(\xi_2) \r|\le \frac{1}{2T},\\
     &\l|\xi_1-\frac{\eta_1}{2} \r|^2 >4 \frac{N_1}{\lambda}, \xi_3=\frac{k}{\lambda} \Big\},
    \end{split}
\end{equation*}
$$ k_0(\xi_2)=\min\{k\in \bZ:   C_k\p(\xi_2) \bigcap \{|(\xi_1,\xi_3)|\sim N_2\}\ne \varnothing , |k|\lesssim \lambda N_2       \},$$
and
$$ k_1(\xi_2)=\max\{k\in \bZ:   C_k\p(\xi_2) \bigcap \{|(\xi_1,\xi_3)|\sim N_2\}\ne \varnothing , |k|\lesssim \lambda N_2       \}.$$
Note that we now have $\mu_{k_0(\xi_2)}\p(\xi_2)> 3 \frac{N_1}{\lambda} $, and for $k\in [k_0(\xi_2),k_1(\xi_2)]$, there holds
$$\mu_{k+1}\p(\xi_2)-  \mu_{k}\p(\xi_2)\sim \frac{N_1}{\lambda}\gg \frac{1}{T},$$
and
$$ |C_k\p(\xi_2)|\lesssim (\lambda T)^{-1}(\mu_{k}\p(\xi_2))^{-\frac12}.$$
If $k_0(\xi_2)=k_1(\xi_2)$,  then we have
\begin{equation}\label{eq:estimate when k0 is k1'}
\sum_{k=k_0(\xi_2)}^{k_1(\xi_2)} |C_k\p(\xi_2| \lesssim (\lambda T)^{-1} (\mu_{k_0(\xi_2)}\p(\xi_2))^{-\frac12} \lesssim \frac{1}{(\lambda N_1)^{\frac12} T}.
\end{equation}
Now we assume $k_1(\xi_2)>k_0(\xi_2)$.   Then for $i=0,1,$ there holds
$$ \l[\l(\mu_{k_i(\xi_2)}\p(\xi_2)-\frac{1}{2T}\r)^{\frac12}, \l(\mu_{k_i(\xi_2)}\p(\xi_2)+\frac{1}{2T}\r)^{\frac12} \r] \bigcap \l\{ \l|\xi_1- \frac{\eta_1}{2}\r|: |\xi_1| \lesssim N_2   \r\} \ne \varnothing,  $$
thus
$$  \l(\mu_{k_1(\xi_2)}\p(\xi_2)-\frac{1}{2T}\r)^{\frac12}-\l(\mu_{k_0(\xi_2)}\p(\xi_2)+\frac1{2T}\r)^{\frac12} \lesssim N_2,$$
so
\begin{equation}\label{eq:estimate k1-k0-2}
\mu_{k_1(\xi_2)}\p(\xi_2) -\mu_{k_0(\xi_2)}\p(\xi_2) \lesssim N_2 (\mu_{k_1(\xi_2)}\p(\xi_2))^{\frac12}.
\end{equation}
When $\mu_{k_0(\xi_2)}\p(\xi_2)\lesssim \mu_{k_1(\xi_2)}\p(\xi_2) -\mu_{k_0(\xi_2)}\p(\xi_2)$, then we get
$$ k_1(\xi_2)-k_0(\xi_2) \lesssim \frac{\lambda N_2^2}{N_1},$$
then
\begin{align*}
\sum_{k=k_0(\xi_2)}^{k_1(\xi_2)} |C_k\p(\xi_2| & \lesssim  \sum_{k=k_0(\xi_2)}^{k_1(\xi_2)} (\lambda T)^{-1}(\mu_{k}\p(\xi_2))^{-\frac12} \\
& \lesssim \sum_{k=k_0(\xi_2)}^{k_1(\xi_2)} (\lambda T)^{-1}  (\mu_{k_0(\xi_2)}\p(\xi_2)+ (k-k_0(\xi_2)) \frac{N_1}{\lambda}  )^{-\frac12} \\
&\lesssim \frac{N_2}{N_1 T}.
\end{align*}
When $\mu_{k_0(\xi_2)}\p(\xi_2)\gg \mu_{k_1(\xi_2)}\p(\xi_2) -\mu_{k_0(\xi_2)}\p(\xi_2)$, by (\ref{eq:estimate k1-k0-2}), one has
$$ \mu_{k_1(\xi_2)}\p(\xi_2) -\mu_{k_0(\xi_2)}\p(\xi_2) \lesssim N_2 (\mu_{k_0(\xi_2)}\p(\xi_2))^{\frac12},$$
so
$$ k_1(\xi_2) -k_0(\xi_2)\lesssim  \frac{\lambda N_2}{N_1} (\mu_{k_0(\xi_2)}\p(\xi_2))^{\frac12} ,$$
thus
\begin{align*}
\sum_{k=k_0(\xi_2)}^{k_1(\xi_2)} |C_k\p(\xi_2| & \lesssim  \sum_{k=k_0(\xi_2)}^{k_1(\xi_2)} (\lambda T)^{-1}(\mu_{k_0(\xi_2)}\p(\xi_2))^{-\frac12} \\
&\sim \l(k_1(\xi_2) -k_0(\xi_2)\r) (\mu_{k_0(\xi_2)}\p(\xi_2))^{-\frac12} \\
&\lesssim  \frac{N_2}{N_1 T}.
\end{align*}

Anyway, we have
$$ \sum_{k=k_0(\xi_2)}^{k_1(\xi_2)} |C_k\p(\xi_2|\lesssim \frac{1}{(\lambda N_1)^{\frac12} T}+ \frac{N_2}{N_1 T}.$$
So there holds
\begin{align*}
&\sum_{k\in \bZ:\mu_k> 10\frac{N_1}{\lambda}}\l| C_k \bigcap \l\{\xi\in \bR^2 \times \bZ_{1/\lambda}:  \l|\xi_1-\frac{\eta_1}{2} \r|^2 >4 \frac{N_1}{\lambda}      \r\}\bigcap  \{|\xi|\sim N_2\} \r| \\
\lesssim & \sum_{k \in \bZ} \int_{|\xi_2|\lesssim N_2} \l| C_k\p(\xi_2) \bigcap  \{|(\xi_1,\xi_3)|\sim N_2\} \r| \dif \xi_2\\
\lesssim & \int_{|\xi_2|\lesssim N_2} \sum_{k=k_0(\xi_2)}^{k_1(\xi_2)} |C_k\p(\xi_2)| \dif \xi_2 \lesssim \frac{N_2}{(\lambda N_1)^{\frac12}T} + \frac{N_2^2}{N_1T} \lesssim \frac{1}{\lambda T}+ \frac{N_2^2}{N_1 T}.
\end{align*}

{\bf Case 2(b).} $m\ge 2, d=m+n\ge 4$. Fix $\xi_i$ with $|\xi_i|\lesssim N_2$, $i=3,\cdots,d-1$, then we may consider $(\xi_1, \xi_2,\xi_d)$ similarly as in Case 2(a).

The proof of \eqref{eq: bilinear2} is now complete.
\hfill$\Box$\vspace{2ex}

\section{Proof of Theorem \ref{thm: nonlinearmain} and Theorem \ref{thm: nonlineark}}
We will focus on the proof of Theorem \ref{thm: nonlinearmain}, and Theorem \ref{thm: nonlineark} follows in a similar way, and we will summarize its proof in the last subsection.
\subsection{Setting up and a recap of I-methods}
It is by now standard to apply I-method, \cite{I-team2}, to obtain low regularity $H^{s}$ GWP for \eqref{eq: cubicnls}. We will mainly follow the presentation\footnote{Strictly speaking, we are working on $\mathbb{R}\times \mathbb{T}$ rather than $\mathbb{T}^{2}$, but the schemes are same. We also have slightly better bilinear estimates, but I-method does not capture $\epsilon$ improvement (or loss). Our better range of $s$ comes from improve the computation in \cite{DPST}.} in \cite{DPST}, but we will improve the computation there to obtain a better range of $s$. We will focus on the different part of the computation, but only sketch the part which are similar.\\

Recall \eqref{eq: cubicnls} conserves its energy,
\begin{equation}
	E(u):=\int \frac{1}{2}|\nabla u|^{2}+\frac{1}{4}|u|^{4}.
\end{equation}
Similarly, one also defines the energy $E_{\lambda}(u)$.

The I-method, \cite{I-team2}, also referred as almost conservation law, is to study modified energy $E(I_{N}u)$, where $I_{N}=m_{N}(D)$ is a multiplier operator, where, $m_{N}(\xi)=\psi(\xi/N)$ and $\psi$ is a smooth function with
\begin{equation}
\psi(\xi)=
\begin{cases}
1, \text{ if } |\xi|\leq 1,\\
|\xi|^{s-1}, \text{ if } |\xi|\geq 2.
\end{cases}
\end{equation}

We note that $I_{N}$ is a smoothed version of $P_{\leq N}$.

It is easy to see $m_{N}$ sends $H^{s}$ to $H^{1}$, thus $E(I_{N}u)$ is well defined for $H^{s}$ functions.

In implementation I-method, one will indeed study \eqref{eq: cubicnls} on $\mathbb{R}\times \mathbb{T}_{\lambda}$,
\begin{equation}\label{eq: rtscale}
\begin{cases}
iu_{t}+\Delta u=|u|^{2}u,\\
u(x,y,0)=f(x,y)\in H^{s}(\mathbb{R}\times \mathbb{T}_{\lambda}).
\end{cases}
\end{equation}

The following proposition has been proved in \cite{DPST}, and can be line by line extended to the $\mathbb{R}\times \mathbb{T}$ case, (for convenience, we only state the $\mathbb{R}\times \mathbb{T}$ version),

\begin{prop}[\cite{DPST}]\label{prop: ac1}
Let $s>\frac{1}{2}$, $1\leq\lambda\lesssim N, 0<t\leq 1$, let $u$ solve \eqref{eq: rtscale}, $\|Iu\|_{X^{1,\frac{1}{2}+}}\lesssim 1$, then
\begin{equation}
E_{\lambda}(I_{N}u(t))-E_{\lambda}(I_{N}u(0))\lesssim \frac{1}{N^{1-}}.	
\end{equation}

\end{prop}
By using standard local theory and  choose $\lambda\sim N^{\frac{1-s}{s}}$, one obtains $H^{s}$ ($s>\frac{2}{3}$) GWP for \eqref{eq: cubicnls} by iterating Proposition \ref{prop: ac1} $N^{1-}$ times, we refer to \cite{DPST,I-team2} for more details.

In current article, we will prove

\begin{prop}\label{prop: ac2}
Let $s>\frac{1}{2}$, $1\leq\lambda\sim N^{\alpha}\lesssim N, 0<t\leq 1$, let $u$ solve \eqref{eq: rtscale}, $\|Iu\|_{X^{1,\frac{1}{2}+}}\lesssim 1$, then
\begin{equation}\label{eq: imethod}
E_{\lambda}(I_{N}u(t))-E_{\lambda}(I_{N}u(0))\lesssim \frac{1}{N^{1+\frac{\alpha}{2}-}}.	
\end{equation}
\end{prop}

By taking $\alpha=\alpha(s)=\frac{1-s}{s}$, i.e.  $\lambda\sim N^{\frac{1-s}{s}}$ as in \cite{DPST}, one obtains $H^{s}$ $(s>\frac{3}{5})$ GWP. See Step 3 in the proof of Theorem 1.2 in \cite{DPST} for more details on how to obtain such an index from Proposition \ref{prop: ac2}. The result can be line by line extended to the $\mathbb{T}^2$ case.

It remains to prove Proposition \ref{prop: ac2}. Since $N$, and $\lambda$ will be fixed throughout, we short $I_{N}, m_{N}$ as $I,M$.

By plug in equation \eqref{eq: rtscale}, one may follow the arguments in \cite{I-team2,DPST}, to derive that $E(Iu(t))-E(Iu(0))$ is bounded by the following two terms,

\begin{equation}
\begin{aligned}
& Tr_1:=\int_0^t \int_{\Gamma_{4}}\left(1-\frac{m\left(\xi_2+\xi_3+\xi_{4}\right)}{m\left(\xi_2\right) m\left(\xi_3\right) m\left(\xi_{4}\right)}\right) \widehat{\Delta \overline{I u}}\left(\xi_1\right) \widehat{I u}\left(\xi_2\right)  \widehat{\overline{I u}}\left(\xi_3\right) \widehat{I u}\left(\xi_{4}\right), \\
& Tr_2:=\int_0^t \int_{\Gamma_{4}}\left(1-\frac{m\left(\xi_2+\xi_3+\cdots+\xi_{4}\right)}{m\left(\xi_2\right) m\left(\xi_3\right)  m\left(\xi_{4}\right)}\right) \widehat{I|u|^{2} \overline{u}}\left(\xi_1\right) \widehat{I u}\left(\xi_2\right) \widehat{\overline{I u}}\left(\xi_3\right) \widehat{I u}\left(\xi_{4}\right),
\end{aligned}
\end{equation}
where $\Gamma_{4}:=\left\{\left(\xi_1, \xi_2, \dots, \xi_{4}\right) \in (\mathbb R\times \mathbb{Z}_{\frac{1}{\lambda}})^{4}: \xi_1+\xi_2+\dots+\xi_{4}=0\right\}$.

It has been proved in \cite{DPST} that
\begin{equation}
|Tr_{2}|\lesssim \frac{1}{N^{2-}}.
\end{equation}
 Thus it remains to prove, given $\|Iu\|_{X^{1,\frac{1}{2}+}}\lesssim 1$
\begin{equation}\label{eq: acmain}
|Tr_{1}|\lesssim \frac{1}{N^{1+\frac{\alpha}{2}-}}.
\end{equation}
As usual, one may localize each term in $Tr_{1}$ at some dyadic frequency $N_{i}$, then sums up. Parallel to (4.37) in \cite{DPST}, in our case, we need to prove
\begin{equation}\label{eq: acworkmain}
    \begin{aligned}
& \left|\int_0^t \int_{\Gamma_{4}}\left(1-\frac{m\left(\xi_2+\xi_3+\xi_{4}\right)}{m\left(\xi_2\right) m\left(\xi_3\right) m\left(\xi_{4}\right)}\right) \widehat{\overline{\phi_1}}\left(\xi_1\right) \widehat{\phi_2}\left(\xi_2\right) \widehat{\overline{\phi_3}}\left(\xi_3\right) \widehat{\phi_{4}}\left(\xi_{4}\right)\right| \\
\lesssim & \frac{1}{N^{\frac{\alpha(s)}{2}+1-}}\max\{\left(N_1 N_2 N_3N_{4}\right)^{0-}, (\frac{N_{1}}{N_{2}}N_{3}N_{4})^{0-}\}\left\|\phi_1\right\|_{X^{-1,1 / 2+}} \prod_{i=2}^{4}\left\|\phi_i\right\|_{X^{1,1 / 2+}},
\end{aligned}
\end{equation}
where $\phi_{i}$ is localised at frequency $\langle \xi_i\rangle \sim N_{i}$, where $N_{i}\geq 1$ is a dyadic number. \\

We note that in \eqref{eq: acworkmain}, when $(N_{1}N_{2}N_{3}N_{4})^{0-}$ dominates, one sums via Minkowski, and when $(\frac{N_{1}}{N_{2}}N_{3}N_{4})^{0-}$ dominates, one sums $N_{3},N_{4}$ via Minkowski, and sums $N_{1},N_{2}$ by Cauchy, (the main part is $N_{1}\sim N_{2}$.)

One may assume $N_{2}\geq N_{3}\geq N_{4}$ by symmetry, (the complex conjugacy will not play a role in the rest part of the  proof).\\

We will still use $Tr_{1}$ to denote the integral on the LHS of \eqref{eq: acworkmain}.

Estimate \eqref{eq: acworkmain} will be analyzed in the Subsection below.
\subsection{Proof of main estimate \eqref{eq: acworkmain}}
One needs to discuss several different cases.
\subsubsection{Case 1. $N\gg N_{2}$}
In this case, by the definition of $m$, one sees the LHS of \eqref{eq: acworkmain} is zero, and \eqref{eq: acmain} holds trivially.
\subsubsection{Case 2. $N_{2}\gtrsim N\gg N_{3}\geq N_{4}$}
Following \cite{DPST}, observing in this case $N_{1}\sim N_{2}$, and using $|1-\frac{m(\xi_{2}+\xi_{2}+\xi_{3})}{m(\xi_{2})}|\lesssim \frac{N_{3}}{N_{2}}$, one obtains
\begin{equation}\label{eq: c2estimate}
|Tr_{1}|\lesssim \frac{N_{3}}{N_{2}}\|\phi_{1}\phi_{3}\|_{L_{t,x,y}^{2}}\|\phi_{2}\phi_{4}\|_{L_{t,x,y}^{2}},
\end{equation}
i.e. (4.40) in \cite{DPST}, and they further estimate this via
\begin{equation}
|Tr_{1}|\lesssim \frac{N_{2}^{0-}}{N^{1-}}\|\phi_{1}\|_{X^{-1\frac{1}{2}+}}\Pi_{i=2}^{4}\|\phi_{i}\|_{X^{1,\frac{1}{2}+}}.
\end{equation}
We do slightly better here.

By applying \eqref{eq: xsbbilinear1}, one further estimate \eqref{eq: c2estimate} via
\begin{equation}\label{eq: im1}
\begin{aligned}
	|Tr_{1}|&\lesssim \frac{N_{3}}{N_{2}}\|\phi_{1}\phi_{3}\|_{L_{t,x,y}^{2}}\|\phi_{2}\phi_{4}\|_{L_{t,x,y}^{2}}\\
	&\lesssim \frac{N_{3}}{N_{2}}\l(\frac{1}{\lambda}+\frac{N_{3}}{N_{1}}\r)^{\frac{1}{2}}\l(\frac{1}{\lambda}+\frac{N_{4}}{N_{2}}\r)^{1/2}\frac{N_{1}}{N_{2}N_{3}N_{4}}\|\phi_{1}\|_{X^{-1,\frac{1}{2}+}}\Pi_{i=2}^{4}\|\phi_{i}\|_{X^{1,\frac{1}{2}+}}.
\end{aligned}
\end{equation}

Since $N_{1}\sim N_{2}$, one has
\begin{equation}\label{eq: c2r1}
	\frac{N_{3}}{N_{2}}\l(\frac{1}{\lambda}+\frac{N_{3}}{N_{1}}\r)^{\frac{1}{2}}\l(\frac{1}{\lambda}+\frac{N_{4}}{N_{2}}\r)^{1/2}\frac{N_{1}}{N_{2}N_{3}N_{4}}\lesssim \frac{1}{N_{2}N_{4}}\l(\frac{1}{\lambda}+\frac{N_{3}}{N_{2}}\r)^{\frac{1}{2}}\l(\frac{1}{\lambda}+\frac{N_{4}}{N_{2}}\r)^{1/2}.
\end{equation}
The RHS of the above is monotone decreasing in $N_{4}$, and $1\leq N_{4}\leq N_{3}$. Thus, (note also $\lambda\sim N^{\alpha}\lesssim N\lesssim N_{2}$,) one estimate \eqref{eq: c2r1} by
\begin{equation}
\frac{1}{N_{2}}\l(\frac{1}{\lambda}+\frac{N_{3}}{N_{2}}\r)^{\frac{1}{2}}\frac{1}{\lambda^{\frac{1}{2}}}\lesssim N^{-(1+\frac{\alpha}{2}-)}(N_{1}N_{2}N_{3}N_{4})^{0-}.	
\end{equation}

This gives our desired improved bound.

\subsubsection{Case 3: $N_{2}\geq N_{3}\gtrsim N$}
Following \cite{DPST}, we will need two facts,
\begin{itemize}
\item
\begin{equation}\label{eq: crude1}
\left|1-\frac{m(\xi_{2}+\xi_{3}+\xi_{4})}{m(\xi_{1})m(\xi_{2})m(\xi_{3})}\right|\lesssim \frac{m(\xi_{1})}{m(\xi_{2})m(\xi_{3})m(\xi_{4})}.	
\end{equation}
\item For $\beta\geq \frac{1}{2}-$,
\begin{equation}\label{eq: crude2}
\frac{1}{m(\xi)|\xi|^{\beta}}\lesssim N^{-\beta}.
\end{equation}
\end{itemize}
Similar to \eqref{eq: c2estimate}, one now estimates, via \eqref{eq: crude1},
\begin{equation}\label{eq: c3estimate}
|Tr_{1}|\lesssim \frac{m(N_{1})}{m(N_{2})m(N_{3})m(N_{4})}\|\phi_{1}\phi_{3}\|_{L_{t,x,y}^{2}}\|\phi_{2}\phi_{4}\|_{L_{t,x,y}^{2}}.
\end{equation}
We again apply \eqref{eq: xsbbilinear1}, and obtain
\begin{equation}\label{eq: c3estimate2}
\begin{aligned}
 &\frac{m(N_{1})}{m(N_{2})m(N_{3})m(N_{4})}\|\phi_{1}\phi_{3}\|_{L_{t,x,y}^{2}}\|\phi_{2}\phi_{4}\|_{L_{t,x,y}^{2}}\\
\lesssim &\frac{m(N_{1})}{m(N_{2})m(N_{3})m(N_{4})}\left(\frac{1}{\lambda}+\frac{\min\{N_3,N_1\}}{\max\{N_3,N_1\}}\right)^{\frac{1}{2}} \left(\frac{1}{\lambda}+\frac{N_4}{N_2}\right)^{\frac{1}{2}} \frac{N_1}{N_2 N_3 N_4} \left\|\phi_1\right\|_{X^{-1,1 / 2+}}\prod_{i=2}^{4}\left\|\phi_i\right\|_{X^{1,1 / 2+}}.
\end{aligned}
\end{equation}
In (4.43), (4.44) of \cite{DPST}, the authors estimate  $\left(\frac{1}{\lambda}+\frac{\min\{N_3,N_1\}}{\max\{N_3,N_1\}}\right)^{\frac{1}{2}} \left(\frac{1}{\lambda}+\frac{N_4}{N_2}\right)^{\frac{1}{2}}$ as 1, (up to some extra $\epsilon$ loss since they are in the torus), and then estimate \eqref{eq: c3estimate2} via $\frac{1}{N^{1-}}(N_{1}N_{2}N_{3}N_{4})^{0-}$. Thus, when $\frac{N_{4}}{N_{2}}\leq \frac{1}{\lambda}$, by plug in $\left(\frac{1}{\lambda}+\frac{\min\{N_3,N_1\}}{\max\{N_3,N_1\}}\right)^{\frac{1}{2}} \left(\frac{1}{\lambda}+\frac{N_4}{N_2}\right)^{\frac{1}{2}}\lesssim \lambda^{-\frac{1}{2}}$, and follow \cite{DPST}, the desired estimate will follow.

When $\frac{N_{4}}{N_{2}}>\frac{1}{\lambda}$, one estimate as

\begin{equation}
\begin{aligned}
&\frac{m(N_{1})}{m(N_{2})m(N_{3})m(N_{4})}\left(\frac{1}{\lambda}+\frac{\min\{N_3,N_1\}}{\max\{N_3,N_1\}}\right)^{\frac{1}{2}} \left(\frac{1}{\lambda}+\frac{N_4}{N_2}\right)^{\frac{1}{2}} \frac{N_1}{N_2 N_3 N_4}\\
\lesssim &\frac{m(N_{1})}{m(N_{2})m(N_{3})m(N_{4})}(\frac{N_{4}}{N_{2}})^{\frac{1}{2}}\frac{N_{1}}{N_{2}N_{3}N_{4}}.
\end{aligned}
\end{equation}
 Using \eqref{eq: crude2} and the fact $\frac{m(N_{1})N_{1}^{1-}}{m(N_{2})N_{2}^{-}}\lesssim 1$, one estimates
 \begin{equation}
 \begin{aligned}
 	&\frac{m(N_{1})}{m(N_{2})m(N_{3})m(N_{4})}(\frac{N_{4}}{N_{2}})^{\frac{1}{2}}\frac{N_{1}}{N_{2}N_{3}N_{4}}\\\lesssim
 	&\frac{m(N_{1})N_{1}^{-}}{m(N_{2})N_{2}^{-}}\frac{1}{N_{2}^{1/2}m(N_{3})N_{3}^{1-}}\frac{1}{N_{3}^{0+}}(\frac{N_{1}}{N_2})^{+}\frac{1}{N_{4}^{1/2-}m(n_{4})}\frac{1}{N_{4}^{+}}\\
 	\lesssim & N^{\frac{3}{2}-}\frac{1}{N_{3}^{0+}}(\frac{N_{1}}{N_2})^{+}\frac{1}{N_{4}^{+}},
 \end{aligned}
 \end{equation}
 which gives the desired bound.

 \subsection{The general case, proof of Theorem \ref{thm: nonlineark}}
 
 The general case follows in a similar way, and this time one needs to prove

 \begin{equation}\label{eq: acworkmain2}
    \begin{aligned}
& \left|\int_0^t \int_{\Gamma_{2k+2}}\left(1-\frac{m\left(\xi_2+\xi_3\cdots\xi_{2k+2}\right)}{m\left(\xi_2\right) m\left(\xi_3\right)\cdots m\left(\xi_{2k+2}\right)}\right) \widehat{\overline{\phi_1}}\left(\xi_1\right) \widehat{\phi_2}\left(\xi_2\right) \widehat{\overline{\phi_3}}\left(\xi_3\right) \cdots \widehat{\phi_{2k+2}}\left(\xi_{2k+2}\right)\right| \\
& \lesssim \frac{1}{N^{\frac{\alpha(s)}{2}+1-}}\left(N_1 N_2 N_3\cdots N_{2k+2}\right)^{0-}\left\|\phi_1\right\|_{X^{-1,1 / 2+}} \prod_{i=2}^{2k+2}\left\|\phi_i\right\|_{X^{1,1 / 2+}},
\end{aligned}
\end{equation}
where $\phi_{i}$ is localised at frequency $\langle \xi_i \rangle \sim N_{i}$, where $N_{i}\geq 1$ is a dyadic number, and $N_{2}\geq N_{3}\cdots \geq N_{2k+2}$, and all the $\phi_{i}$, $i\geq 4$ will be estimated via $L_{t,x}^{\infty}$, which is bounded by its $X^{1,\frac{1}{2}+}$ norm by Sobolev embedding.

{Note that in both cases, \eqref{eq: cubicnls}, and \eqref{eq: knls}, we get the same estimates on the rescaled waveguide \eqref{eq: imethod} formally with $\alpha(s):=\frac{s-1}{1-s-1 / k}$ , just regarding \eqref{eq: cubicnls} as a special case of \eqref{eq: knls} by taking $k=1$. In addition, \eqref{eq: imethod} holds for \eqref{eq: knls} at the range of $s> 1-\frac{1}{2k}$. However, after undoing the scaling, we use \eqref{eq: imethod} to iterate constructing the solution of \eqref{eq: knls} on $[0,T]$,  and get $\|u(T)\|_{H^s} \lesssim C_{N,\lambda}$ for $T\ll N^{\frac{\frac{5}{2}(1-s)-\frac{1}{k}}{1-s-\frac{1}{k}}}$.  This  leads to the range of low regularity $s>1-\frac{2}{5k}$ GWP for \eqref{eq: knls}. }\\

One subtle thing here is now the crude bound \eqref{eq: crude1} will be replaced via
\begin{equation}
|1-\frac{m(\xi_{2}+\cdots m(\xi_{2k+2}))}{m(\xi_{2})\cdots m(\xi(2k+2))}|\lesssim \frac{m(\xi_{1})}{m(\xi_{2})\cdots m(\xi_{j_0})}	.
\end{equation}
Unless there exists $4\leq j_0\leq 2k+2$, such that $N_{j_0} \gtrsim N\gg N_{j_0+1}$(with a little abuse, when $j_0=2k+2$, $N_{2k+2} \gtrsim N\gg N_{2k+2}$ represents the case $N_2\geq N_3\geq\cdots\geq N_{2k+2}\gtrsim N$ ).  This part will need extra treatment.

% Similarly, we will use the following bound on the multiplier,
%$$
%\left|1-\frac{m\left(\xi_2+\xi_3+\cdots+\xi_{2k+2}\right)}{m\left(\xi_2\right) m\left(\xi_3\right) \cdots m\left(\xi_{2k+2}\right)}\right| = \left|1-\frac{m\left(\xi_2+\xi_3+\cdots+\xi_{2k+2}\right)}{m\left(\xi_2\right) m\left(\xi_3\right) \cdots m\left(\xi_{j_0}\right)}\right|\lesssim  \frac{m\left(\xi_1\right)}{m\left(\xi_2\right) m\left(\xi_3\right) \cdots m\left(\xi_{j_0}\right)}.
%$$

Through a rough estimate $\left(\frac{1}{\lambda}+\frac{N_3}{N_1}\right)^{\frac{1}{2}} \left(\frac{1}{\lambda}+\frac{N_4}{N_2}\right)^{\frac{1}{2}} \lesssim 1$, one now estimate via
$$
\begin{aligned}
\left|\operatorname{Tr}_1\right| \lesssim & \frac{m\left(N_1\right)}{\prod\limits_{j=2}^{j_0} m\left(N_j\right) }
%\lambda^{0+} \frac{\left(N_1\right)^{1+}}{\left(N_2 N_3 N_4\right)^{1-}}
\left\| \phi_{1}\phi_{3}\right\|_{L^2_{t,x}}\left\|\phi_{2}\phi_{4}\right\|_{L^2_{t,x}}\prod_{i=5}^{2k+2}\|\phi_i\|_{L^\infty_{t,x}} \\
%\lesssim&  \frac{m\left(N_1\right)}{\prod\limits_{j=2}^{j_0} m\left(N_j\right) } \left(\frac{1}{\lambda}+\frac{N_3}{N_1}\right)^{\frac{1}{2}} \left(\frac{1}{\lambda}+\frac{N_4}{N_2}\right)^{\frac{1}{2}} \frac{N_1}{N_2 N_3 N_4} \left\|\phi_1\right\|_{X^{-1,1 / 2+}}\prod_{i=2}^{2k+2}\left\|\phi_i\right\|_{X^{1,1 / 2+}}\\
% \lesssim &  \frac{m(N_1)N_1}{m(N_2)N_2}\frac{1}{N_3 N_4\prod^{j_0}_{j=3}m\left(N_j\right) }\left\|\phi_1\right\|_{X^{-1,1 / 2+}}\prod_{i=2}^{2k+2}\left\|\phi_i\right\|_{X^{1,1 / 2+}}\\
 \lesssim &  \frac{m(N_1)N_1^{1-}}{m(N_2)N_2^{1-}}\frac{1}{ m(N_3)N_3^{1-}}\frac{1}{ N_4\prod^{j_0}_{j=4}m\left(N_j\right) } \frac{N_1^{0+}}{N_2^{0+}} \frac{1}{N_3^{0+}}\left\|\phi_1\right\|_{X^{-1,1 / 2+}}\prod_{i=2}^{2k+2}\left\|\phi_i\right\|_{X^{1,1 / 2+}}\\
% \lesssim & \frac{1}{N^{1-}} \frac{1}{\prod^{j_0}_{j=4}m\left(N_j\right) N_4} \frac{N_1^{0+}}{N_2^{0+}} \frac{1}{N_3^{0+}}\left\|\phi_1\right\|_{X^{-1,1 / 2+}}\prod_{i=2}^{2k+2}\left\|\phi_i\right\|_{X^{1,1 / 2+}}\\
 % \lesssim & \frac{1}{N^{1-}} \frac{1}{\prod^{j_0}_{j=4}\left(m\left(N_j\right) N_j^{\frac{1}{j_0-3}}\right)} \frac{N_1^{0+}}{N_2^{0+}} \frac{1}{N_3^{0+}}\left\|\phi_1\right\|_{X^{-1,1 / 2+}}\prod_{i=2}^{2k+2}\left\|\phi_i\right\|_{X^{1,1 / 2+}}\\
  \lesssim & \frac{1}{N^{2-}} \frac{N_1^{0+}}{N_2^{0+}} \frac{1}{N_3^{0+}}\left\|\phi_1\right\|_{X^{-1,1 / 2+}}\prod_{i=2}^{2k+2}\left\|\phi_i\right\|_{X^{1,1 / 2+}},
\end{aligned}
$$provided that $s> 1-\frac{1}{2k-1}$. 
The factor $\frac{N_1^{0+}}{N_2^{0+}} \frac{1}{N_3^{0+}}$ allows us to directly sum in $N_3, N_4, \ldots, N_{2 k+2}$, and sum in $N_1$ and $N_2$ after applying Cauchy-Schwarz to those factors.

\appendix
\section{On sharpness of estimates \eqref{eq: bilinear1}, \eqref{eq: bilinear2}}\label{app: sharp}
\subsection{Sharpness of the local-in-time bilinear estimate}
We first take some examples to show the sharpness of the local-in-time estimate
\begin{equation} \label{appendix-eq-local-in-time bilinear estimate}
\| U_\lambda(t) P_{N_1}f \cdot U_\lambda(t) P_{N_2}g   \|_{L^2(\bR^m \times \bT^n_\lambda\times [0,1])} \lesssim \l( K(\lambda,N_1,N_2)\r)^{\frac12} \|f\|_{L^2} \|g\|_{L^2}, N_1 \ge N_2,
\end{equation}
where
\begin{equation*}
K(\lambda,N_1,N_2):=\left\{\begin{array}{l}
\frac{1}{\lambda} +\frac{N_2}{N_1}, m=n=1, \\
 \frac{N_2^{d-3}}{\lambda} +\frac{N_2^{d-1}}{N_1}, d \geq 3.
\end{array}\right.
\end{equation*}
These examples essentially appeared in \cite{FSWW}. Since we are in the waveguide case rather than pure tori case, we present the examples below for the connivence of readers.\\
%For simplicity, we denote $\xi=(\widetilde{\xi},\xi_d)\in (\bR^m \times \bZ_{1/\lambda}^{n-1} ) \times \bZ_{1/\lambda}$.

We take
$$
U_\lambda(t) P_{N_1}f(x)=\int_{|(\xi_2,\cdots,\xi_d)|\lea N_2, \atop |\xi_1-N_1|\sim N_2} e^{2 \pi i x \cdot \xi-|2 \pi \xi|^2 i t}  ({\rm d} \xi)_\lambda,
$$
and
$$
U_\lambda(t) P_{N_2}g(x)=\int_{|(\xi_2,\cdots,\xi_d)|\lea N_2, \atop |\xi_1|\sim N_2} e^{2 \pi i x \cdot \xi-|2 \pi \xi|^2 i t}  ({\rm d} \xi)_\lambda,
$$
then $|U_\lambda(t) P_{N_2}g(x)|\gea N_2^d$ for $(x,t)\in B:= \{ |x|\lea \frac{1}{N_2}, |t|\lea \frac{1}{N_2^2}   \}$, and  $|U_\lambda(t) P_{N_1}f(x)|\gea N_2^d$ for $(x,t)\in A:=\{ |(x_2,\cdots, x_{d})| \lea \frac{1}{N_2}, |x_1-4\pi N_1 t|\lea \frac{1}{N_2}, |t|\lea \frac{1}{N_2^2}   \}$. Note that $|A\bigcap B|\gea \frac{1}{N_1 N_2^{d+1}} $, so (\ref{appendix-eq-local-in-time bilinear estimate}) reduces to
$$  N_2^{2d-\frac{d+1}{2}} N_1^{-\frac12} \lea (K(\lambda,N_1,N_2))^{\frac12} N_2^d   ,$$
that is
$$ K(\lambda,N_1,N_2)\gea \frac{N_2^{d-1}}{N_1}.$$

When $m=n=1$,  we consider
$$
U_\lambda(t) P_{N_1}f(x)=\int_{|\xi_1|\lea 1, \atop \xi_2=N_1} e^{2 \pi i x \cdot \xi-|2 \pi \xi|^2 i t}  ({\rm d} \xi)_\lambda,
$$
and
$$
U_\lambda(t) P_{N_2}g(x)=\int_{|\xi_1|\lea 1, \atop \xi_2=N_2} e^{2 \pi i x \cdot \xi-|2 \pi \xi|^2 i t}  ({\rm d} \xi)_\lambda,
$$
then $|U_\lambda(t) P_{N_2}g(x)|\gea \frac{1}{\lambda}$ and $|U_\lambda(t) P_{N_1}f(x)|\gea \frac{1}{\lambda}$ for $(x,t)\in \{ |x_1|\lea 1, x_2\in \bT_\lambda, |t|\lea 1 \}$, so (\ref{appendix-eq-local-in-time bilinear estimate}) reduces to
$$    \lambda^{-\frac{3}{2}} \lea  (K(\lambda,N_1,N_2))^{\frac12} \lambda^{-1},           $$
that is
$$ K(\lambda,N_1,N_2)\gea \frac{1}{\lambda}.$$

When $m\ge 1, n\ge 1, d=m+n\ge 3$, we consider
$$
U_\lambda(t) P_{N_1}f(x)=\int_{|(\xi_1,\cdots,\xi_{d-1})|\lea N_2, \atop \xi_d=N_1} e^{2 \pi i x \cdot \xi-|2 \pi \xi|^2 i t}  ({\rm d} \xi)_\lambda,
$$
and
$$
U_\lambda(t) P_{N_2}g(x)=\int_{|(\xi_1,\cdots,\xi_{d-1})|\lea 1, \atop \xi_d=N_2} e^{2 \pi i x \cdot \xi-|2 \pi \xi|^2 i t}  ({\rm d} \xi)_\lambda,
$$
then $|U_\lambda(t) P_{N_2}g(x)|\gea \frac{N_2^{d-1}}{\lambda}$ and $|U_\lambda(t) P_{N_1}f(x)|\gea \frac{N_2^{d-1}}{\lambda}$ for $(x,t)\in \{ |(x_1,\cdots,x_{d-1})|\lea \frac{1}{N_2}, x_d\in \bT_\lambda, |t|\lea \frac{1}{N_2^2} \}$, so (\ref{appendix-eq-local-in-time bilinear estimate}) reduces to
$$    N_2^{2d-2-\frac{d+1}{2}} \lambda^{-\frac{3}{2}} \lea  (K(\lambda,N_1,N_2))^{\frac12} N_2^{d-1} \lambda^{-1},           $$
that is
$$ K(\lambda,N_1,N_2)\gea \frac{N_2^{d-3}}{\lambda}.$$

\subsection{Global-in-time bilinear estimate cannot hold for $m=1, n\ge 1$}
We prove in this subsection, for waveguide with only one $\mathbb{R}$ directions, and $N_{1}\gg N_{2}$, one can construct $f,g\in L^2(\bR \times \bT_\lambda^n)$ so that 
\begin{equation}\label{eq: Global-in-time bilinear estimate}
\|U_{\lambda}(t)P_{N_{1}}f U_{\lambda}(t)P_{N_{2}}g\|_{L^2(\mathbb{R}\times \mathbb{T}_\lambda^n \times [0,\infty])}=\infty, N_{1}\gg N_{2}.	
\end{equation}

Note that, when there is no periodic part, in 1d Euclidean case, one does have global bilinear estimates
$$
\|e^{it\Delta_\bR} P_{N_{1}}f e^{i\Delta_\bR}P_{N_{2}}g\|_{L_{t,x}^{2}(\mathbb{R}\times \mathbb{R})}\lesssim (\frac{1}{N_{1}})^{1/2}\|f\|_{L_{x}^{2}}\|g\|_{L_{x}^{2}}.	
$$

On the other hand, it is not hard to see, in general there cannot hold 
\begin{equation}
\|e^{it\Delta_\bR} P_{\leq 1}f e^{i\Delta_\bR}P_{\leq 1}g\|_{L_{t,x}^{2}(\mathbb{R}\times \mathbb{R})}<\infty
\end{equation}
for all $f,g$ Schwarz.

From this perspective, in the case $\bR \times \bT_\lambda^n$, the frequencies of $f$ and $g$ may concentrate on  $\bZ_\lambda^n$ part, thus the assumption $N_{1}\gg N_{2}$ can not ensure frequencies of $f$ and $g$ on $\bR$ be separated.

For simplicity, we consider the case $\mathbb{R}\times \mathbb{T}_{\lambda}$, and the other cases are similar.

Let

\begin{equation}
\begin{aligned}
&U_\lambda(t) P_{N_1}f(x)=\int_{ \xi_2=N_1} e^{2 \pi i x \cdot \xi-|2 \pi \xi|^2 i t} \widehat{\varphi}(\xi_1) ({\rm d} \xi)_\lambda,\\
&U_\lambda(t) P_{N_2}g(x)=\int_{ \xi_2=N_2} e^{2 \pi i x \cdot \xi-|2 \pi \xi|^2 i t}  \widehat{\varphi}(\xi_1) ({\rm d} \xi)_\lambda,
\end{aligned}
\end{equation}
where   $\varphi= P_{<1}(e^{-y^2})$, one observes 
\begin{equation}
\|U_{\lambda}(t)P_{N_{1}}f U_{\lambda}(t)P_{N_{2}}g\|_{L^2(\mathbb{R}\times \mathbb{T}_\lambda^n \times [0,\infty])}\sim_{\lambda}\|e^{it \Delta_\bR} \varphi \|_{L^4(\bR \times [0,\infty])}^{2}.
\end{equation}

It suffices to show that
\begin{equation}
 \|e^{it \Delta_\bR} \varphi \|_{L^4(\bR \times [0,\infty])}=\infty.
 \end{equation}
It follows from that $ |e^{it \Delta_\bR} \varphi(x)|\gea t^{-\frac12}  $ when $|x|\le \frac{t}{1000}$ for $t\gg 1$.

\section{A sketch of LWP for \eqref{eq: knls}}
We fix $k$ and also fix $s=s_{c}=1-\frac{1}{k}$.

By \cite{TT}, one has 
\begin{equation}
\|e^{it\Delta}P_{N}f\|_{L_{t,x,y}^{4}(\mathbb{R}\times \mathbb{T}\times [0,1])}\leq \|P_{N}f\|_{L^{2}}.
\end{equation}

Note that by Sobolev embedding
\begin{equation}
\|e^{it\Delta}P_{N}f\|_{L_{t,x,y}^{\infty}(\mathbb{R}\times \mathbb{T}\times [0,1])}\leq N\|P_{N}f\|_{L^{2}}.
\end{equation}

By interpolation, one has
 \begin{equation}
\|e^{it\Delta}P_{N}f\|_{L_{t,x,y}^{4k}(\mathbb{R}\times \mathbb{T}\times [0,1])}\leq N^{s}\|P_{N}f\|_{L^{2}}.
\end{equation}
Thus, one has 
\begin{equation}
	\|e^{it\Delta}f\|_{L_{t,x,y}^{4k}(\mathbb{R}\times \mathbb{T}\times [0,1])}\leq \|\langle \nabla \rangle^{s}f\|_{L^{2}}.
\end{equation}

Based on above estimate, one may prove LWP of \eqref{eq: knls} under norms $\|\langle \nabla \rangle^{s}f\|_{L_{t,x,y}^{4}}$, $\|f\|_{L_{t}^{\infty}H^{s}}$ and $\|f\|_{L^{4k}_{t,x,y}}$ by Picard iterations.

Let $\|f\|_{X}:=\|\langle \nabla \rangle^{s}f\|_{L_{t,x,y}^{4}}+\|f\|_{L_{t}^{\infty}H^{s}}+\|f\|_{L^{4k}_{t,x,y}}.$

We list the needed estimates below, and leave the rest to interested readers.
\begin{enumerate}
	\item For linear solutions, one applies
	\begin{equation}
		\|f\|_{X}\lesssim \|f\|_{H^{s}}.
		\end{equation}
    \item For Duhamel part, one applies, (dual estimates),
    \begin{equation}
    \|\int_{0}^{t}e^{i(t-s)\Delta}F\|_{X}\lesssim \|\langle \nabla \rangle^{s} F\|_{L_{t}^{4/3}L_{x,y}^{4/3}}.	
    \end{equation}
\item And, finally, nonlinear estimates
\begin{equation}
	\|\langle \nabla \rangle^{s}|u|^{2k}u\|_{L_{t}^{4/3}L_{x,y}^{4/3}}\lesssim \|\langle \nabla \rangle^{s}u\|_{L_{t,x,y}^{4}}\|u\|_{L^{4k}_{t,x,y}}^{2k}.
\end{equation}

\end{enumerate}

\bibliographystyle{amsplain}
\bibliographystyle{plain}
\bibliography{BG}
\end{document}